\ifx\LocalElevenPtArticleFormatLoaded\undefined
  \documentclass[a4paper,11pt]{article}     % `article' for A4 paper
  \usepackage{array}                        % array extensions
  \usepackage{theorem}                      % theorem extensions
  \usepackage{amsmath,amscd,amssymb}        % AMS-LaTeX with CD and symbols
  \usepackage{latexsym}                     % old LaTeX symbols
  \usepackage[german]{babel}                % German only
  \usepackage{graphicx}                     % graphic extension for picture
  \usepackage{epsfig}
\fi

\usepackage{theorem}
\newtheorem{theorem}{Satz}[section]
\newtheorem{lemma}[theorem]{Lemma}

\newtheorem{cor}[theorem]{Korollar}
\theorembodyfont{\rm}

\newtheorem{definition}{Definition.} 
      
\newtheorem{uremark}{Bemerkung.}

\newcommand{\CC}{{\mathbb{C}}}

\newcommand{\HH}{{\mathbb{H}}}
\newcommand{\LL}{{\mathbb{L}}}

\newcommand{\PP}{{\mathbb{P}}}
\newcommand{\QQ}{{\mathbb{Q}}}
\newcommand{\ZZ}{{\mathbb{Z}}}

\newcommand{\bbQ}{{\mathbb{Q}}}
\newcommand{\RR}{{\mathbb{R}}}

\newcommand{\bbEins}{{\mathbf{1}}}

\newcommand{\cO}{{\cal O}}

\newcommand{\cM}{{\cal M}}
\newcommand{\cN}{{\cal N}}

\newcommand{\cD}{{\cal D}}

\newcommand{\on}[1]{\operatorname{#1}}
\newcommand{\Gpq}{{\tilde{\Gamma}_{1,p}(q)}}
\newcommand{\Gp}{{\tilde{\Gamma}^{\circ}_{1,p}}}
\newcommand{\Gt}{{\tilde{\Gamma}^{\circ}_{1,t}}}
\newcommand{\Gtn}{{\tilde{\Gamma}_{1,t}(n)}}
\newcommand{\Gtlev}{{\tilde{\Gamma}_{1,t}}}
\newcommand{\Gpp}{{\tilde{\Gamma}_{1,p}(p)}}
\newcommand{\pmat}[1]{\begin{pmatrix}#1\end{pmatrix}}
\newcommand{\gh}[1]{\mathfrak{H}\left(\left[#1\right]_{\cO(t)}\right)}
\newcommand{\gl}[1]{\mathfrak{L}\left(\left[#1\right]_{\cM(r)}\right)}
\renewcommand{\Sp}{\on{Sp}}
\newcommand{\Sym}{\on{Sym}}

\newcommand{\Gr}{\on{Gr}}
\newcommand{\prim}{\on{prim}}
\newcommand{\diag}{\on{diag}}
\newcommand{\ggT}{\on{ggT}}
\newcommand{\GL}{\on{GL}}
\newcommand{\SL}{\on{SL}}

\newenvironment{Proof*}[1]{\begin{ProofwCaption}{{#1}}}{\end{ProofwCaption}}
\newenvironment{ProofwCaption}[1]%
  {\addvspace\theorempreskipamount \noindent{\it #1.}\rm}%
  {\qed \par \addvspace\theorempostskipamount}
\newcommand{\qedsymbol}{\mbox{$\Box$}}
\newcommand{\qed}{\quad\qedsymbol}
\begin{document}
\title{Das Titsgeb\"aude von Siegelschen Modulgruppen vom Geschlecht~2}
\author{M.~Friedland und G.K.~Sankaran}
\date{}
\maketitle

\begin{quotation}
\begin{center}
{\bf Abstract}
\end{center}

  \noindent We describe the Tits buildings of the Siegel modular
  groups $\Gt$ and $\Gtlev$ for $t$ squarefree. These give the
  configuration of boundary components in the moduli spaces of abelian
  surfaces with polarization of type $(1,t)$ and polarization of type
  $(1,t)$ with canonical level structure. We also do the same
  calculation for the Siegel modular groups $\Gpq$ for $p,q$ prime,
  $p$ odd, corresponding to moduli of abelian surfaces with
  polarization of type $(1,p)$ and full level-$q$ structure. Colour
  pictures are available at:
  \begin{center}
    {\tt http://www.bath.ac.uk/\~{}masgks/Buildings}
  \end{center}
\end{quotation}

Wir werden die Titsgeb\"aude gewisser Siegelscher Modulgruppen, die
Mo\-dulr\"aumen polarisierter abelscher Variet\"aten mit Levelstruktur
entsprechen, berechnen. Diese Geb\"aude geben bekanntlich die Struktur
des Randes des Modulraums wieder. Die Berechnungen sind v\"ollig
elementar aber etwas l\"astig. Wir wollen erreichen, da\ss{}
k\"unftige Untersuchungen durch diesen Artikel verk\"urzt werden
k\"onnen. Dar\"uberhinaus sind wohl die Aussagen eleganter als die
Beweise.

Im ersten Teil beschreiben wir kurz die Mo\-dulr\"aume und den
Zusammenhang mit Titsgeb\"auden Symplektischer Gruppen. Dem zweiten
bzw.~dritten Teil sind den Berechnungen des Titsgeb\"audes f\"ur $\Gt$
und $\Gtlev$ ($t$ quadratfrei) bzw.~$\Gpq$, $p,q$ prim, $p>2$,
gewidmet.  F\"ur die Gruppen verwenden wir die Notation
aus~\cite{HKW}. Dabei geh\"ort $\Gt$, bzw. $\Gtlev$ zu den
Mo\-dulr\"aumen abelscher Fl\"achen mit Polarisierung vom Typ $(1,t)$
ohne Levelstruktur bzw.~mit kanonischer Levelstruktur und $\Gpq$ zu
den Mo\-dulr\"aumen mit Polarisierung vom Typ $(1,p)$ mit ganzer
Levelstruktur der Stufe~$q$.

Weitere, farbige Bilder befinden sich auf der Seite:
  \begin{center}
    {\tt http://www.bath.ac.uk/\~{}masgks/Buildings}
  \end{center}

Dieser Artikel beruht im wesentlichen auf der Hannoverschen
Diplomarbeit des ersten Autors, die unter Begleitung von
Prof.~K.~Hulek zustande gekommen ist. Die Verfasser danken dem DAAD
und dem British Council f\"ur finanzielle Unterst\"utzung im Rahmen
des ARC-Projekts 313-ARC-XIII-99/45.

\section {Mo\-dulr\"aume polarisierter abelscher Fl\"achen}

Mo\-dulr\"aume polarisierter abelscher Variet\"aten erh\"alt man als
Quotienten der Siegelschen oberen Halbebene $\HH_g$ nach arithmetischen
Untergruppen $\Gamma$ der symplektischen Gruppe $\Sp(2g,\QQ)$. Der
resultierende Quotientenraum $\Gamma\backslash\HH_g$ ist eine
quasiprojektive Variet\"at mit schlimmstenfalls Quotientensingularit\"aten.
Die Frage nach der Kompaktifizierung dieses Raumes stellt sich in
nat\"urlicher Weise. Eine M\"oglichkeit, dieses Problem anzugehen, besteht
in Mumford's Konzept der toroidalen Kompaktifizierung, wie sie in
\cite{HKW} f\"ur den Fall der $(1,p)$-polarisierten Fl\"achen ($p$ prim)
beschrieben wurde. Dabei wird $\HH_2$ verm\"oge der {\it
Cayley--Abbildung} isomorph auf das dreidimensionale beschr\"ankte Gebiet
$\cD_2 :=\left\{ Z \in \Sym(2,\CC): \bbEins - Z \bar Z > 0 \right\}$
abgebildet. F\"ur den topologischen Abschlu\ss{} definiert man dann
sogenannte {\it rationale Randkomponenten}. Es stellt sich heraus, da\ss{}
diese in 1:1 Beziehung zu isotropen Unterr\"aumen von $\QQ^{4}$ stehen.
Aber mehr noch: Eine Randkomponente $F'$ ist genau dann echt im
Abschlu\ss{} einer Randkomponente $F$ enthalten, wenn der zu $F'$ zugeh\"orige
Unterraum $U_{F'}$ den Unterraum $U_{F}$ zu $F$ echt enth\"alt. Die
Operation von $\Sp(4,\RR)$ auf $ \HH_2$ induziert eine Operation auf
$\cD_2$, die in nat\"urlicher Weise auf den topologischen Abschlu\ss{}
$\bar{\cD}_2$ fortgesetzt werden kann.

Interessieren wir uns f\"ur eine arithmetische Untergruppe $\Gamma \subset
\Sp(4,\QQ)$, m\"ussen wir den Quotienten $\cD_2$ nach $\Gamma$
betrachten. Die rationalen Randkomponenten von $\cD_2$ modulo $\Gamma$
stehen dann in 1:1 Beziehung zu den Bahnen isotroper Unterr\"aume von
$\QQ^4$ modulo $\Gamma$, wobei $\Gamma$ auf $\QQ^4$ durch 
\[
	\gamma: v \mapsto v \cdot \gamma
\]
operiert ($v \in \QQ^4$ sei hier und im weiteren Verlauf als Zeilenvektor notiert).

Die Beschreibung der Bahnen isotroper Unterr\"aume modulo $\Gamma$ und
deren Nachbarschaftsrelationen werden in einem Graphen, dem
Titsgeb\"aude der jeweiligen arithmetischen Untergruppe von
$\Sp(4,\QQ)$, kodiert. Dabei ent\-spricht eine Ecke $e(U\cdot\Gamma)$
dieses Graphen einer Bahn eines nichttrivialen iso\-tropen Unterraums
$U \subset \QQ^4$ bez\"uglich der Operation von $\Gamma$ und je zwei
Ecken $e(U_1\cdot\Gamma)$ und $e(U_2\cdot\Gamma)$ werden genau dann
mit einer Kante verbunden, wenn es ein $\gamma\in\Gamma$ gibt, so
da\ss{} $U_1\cdot\gamma \subset U_2$ oder $U_2 \subset U_1\cdot\gamma$ gilt.
Das Titsgeb\"aude gibt also die Konfiguration der Randkomponenten des
jeweils betrachteten Modulraums wieder.

Es seien $\epsilon_1$ und $\epsilon_2$ Basisvektoren des $\CC^2$ und
$\omega_1, \ldots, \omega_4$ eine Basis eines Gitters $L \subset \CC^2$ von
maximalem Rang und $\omega_i = \omega_{1i}\epsilon_1 +
\omega_{2i}\epsilon_2,~i=1,\ldots,4$. Dann ist $\CC^2/ L$ ein komplexer
Torus, der keineswegs projektiv-algebraisch sein mu\ss{}. Ein komplexer
Torus $X = \CC^2/ L$, der gerade diese Eigenschaft hat, hei\ss{}t {\it
abelsche Fl\"ache}. Dies ist genau dann der Fall, falls es eine nicht
entartete alternierende Bilinearform $\alpha$ gibt, die bez\"uglich der
gew\"ahlten Basis durch eine Matrix $A \in \on{Mat}(4,\ZZ)$ dargestellt
wird, so da\ss{} bez\"uglich der Periodenmatrix von $X$
\[
\Pi := \pmat{\omega_{11} & \omega_{12} & \omega_{13} & \omega_{14}\cr
\omega_{21} & \omega_{22} & \omega_{23} & \omega_{24}}
\]
die {\it Riemannschen Relationen}
\[
\begin{array}{lrcl}
\mathrm{(i)} & \Pi A^{-1}\vphantom{\Pi}^t\Pi & = & 0~\mbox{ und} \cr
\mathrm{(ii)}& i\Pi A^{-1}\vphantom{\Pi}^t\bar{\Pi} &>&0
\end{array}
\]
gelten.

F\"ur eine geeignete Wahl der Basis von $L$ kann man erreichen, da\ss{}
$\alpha$ durch die Matrix
\[
\Lambda = \pmat{0 & E \cr -E & 0}~\mbox{ mit } E:= \diag(e_1,e_2)
\]
dargestellt wird. Dabei sind $e_1,e_2$ eindeutig bestimmte positive
ganze Zahl\-en mit $e_1|e_2$. Das Paar $(e_1,e_2)$ hei\ss{}t {\it Typ}
der Polarisierung. Da die Mo\-dulr\"aume abelscher Fl\"achen mit
$(e_1,e_2)$-- und $(ke_1,ke_2)$--Polarisierung kan\-onisch
isomorph sind, darf man ohne Einschr\"ankung $e_1=1$ annehmen. Um die
kombinatorische Aufgabe im angemessenen Rahmen zu halten, werden wir
im weiteren Verlauf nur Polarisierungen vom Typ $(1,t),t>0$ und $t$
quadratfrei betrachten.

\begin{definition}
  Die Gruppe von linearen Automorphismen auf $L$, die die Form $\Lambda$
  invariant l\"a\ss{}t, also
\[
\Sp(\Lambda,\ZZ):=\left\{\gamma\in\GL(4,\ZZ)\mid
\gamma\Lambda\vphantom{\gamma}^t{\gamma}=\Lambda\right\}
\]
  hei\ss{}t {\it symplektische Gruppe bez\"uglich $\Lambda$}.
\end{definition}
Die Operation von $\Sp(\Lambda,\ZZ)$ auf dem Gitter $L$ induziert eine
Operation auf $\HH_2$, die f\"ur eine Matrix $\pmat{A & B \cr C & D} \in
\Sp(\Lambda,\ZZ)$ ($A,B,C,D$ sind $2 \times 2$ Bl\"ocke) durch
\[
\pmat{A&B \cr C & D} : \tau \longmapsto (A\tau + BE)(C\tau + DE)^{-1}E
\]
gegeben ist. Durch Austeilen dieser Operation wird der Quotient
\[
\Sp(\Lambda,\ZZ)\backslash \HH_2
\]
zu einem Modulraum $(e_1,e_2)$-polarisierter abelscher Fl\"achen.

Mit $L^{\vee}$ bezeichnen wir das {\it duale Gitter} zu $L$ bez\"uglich
$\alpha$, das hei\ss{}t
\[
L^{\vee} = \left\{ y \in L \otimes_{\ZZ} \RR \mid \alpha (x,y) \in \ZZ \mbox{
  f\"ur alle } x \in L \right\}
\]
Der Quotient $L^{\vee} / L$ ist eine endliche Gruppe, die isomorph zu
$(\ZZ_{e_1} \times \ZZ_{e_2})^2$ ist. Diese tr\"agt eine alternierende Form
$\alpha'$, die in den kanonischen Erzeugenden durch die Matrix
\[
        \pmat{0 & E^{-1} \cr -E^{-1} & 0}
\]
bestimmt ist. Eine {\it Levelstruktur vom kanonischen Typ} ist ein
Isomorphismus
\[
\lambda:L^{\vee}/L \stackrel{\sim}{\longrightarrow}(\ZZ_{e_1} \times
\ZZ_{e_2})^2
\]
der die durch $\alpha$ auf $L^{\vee}/L$ induzierte alternierende Form in
die Form $\alpha'$ \"uberf\"uhrt.

Wird nun zus\"atzlich die Erhaltung der Levelstruktur gefordert, so gibt es
wie zuvor eine Korrespondenz zwischen den Punkten von $\HH_2$ und abelschen
Fl\"achen mit Levelstruktur, jedoch verkleinert sich die Gruppe der
Automorphismen zu einer Untergruppe von $\Sp(\Lambda,\ZZ)$.

Setzen wir $\LL = \ZZ^4$, dann ist $\LL^{\vee} = \frac{1}{e_{1}} \ZZ
\oplus \frac{1}{e_{2}} \ZZ \oplus \frac{1}{e_{1}} \ZZ \oplus
\frac{1}{e_{2}} \ZZ$ bez\"uglich $\Lambda$. Die Automorphismen in
$\Sp(\Lambda,\ZZ)$, die die Identit\"at auf $\LL^{\vee}/\LL$
induzieren (und damit die Levelstruktur erhalten), sind gerade die
Elemente $g$ aus $\Sp(\Lambda,\ZZ)$, f\"ur die $vg \equiv v\mod \LL$
f\"ur alle $v \in \LL^{\vee}$ gilt.

Zus\"atzlich betrachten wir das Gitter $\LL_{n}^{\vee} = \frac{1}{n}
\ZZ \oplus \frac{1}{n} \ZZ \oplus \frac{1}{n} \ZZ \oplus \frac{1}{n}
\ZZ$, das duale Gitter von $\LL$ bez\"uglich $nJ:=n\pmat{0&\bbEins_2\\
  -\bbEins_2&0}.$

\begin{definition}
Wir definieren folgende Matrixgruppen:
\[
        \begin{array}{ccl}
                \Gt &=& \Sp(\Lambda,\ZZ)\\
                \Gtlev &=& \left\{g \in \Gt\mid vg \equiv v \mod \LL~ \mbox{ f\"ur
                        alle } v \in \LL^{\vee}\right\} \\
                \Gtn &=& \left\{g \in \Gt \mid vg \equiv v \mod \LL~\mbox{ f\"ur
                        alle } v \in \LL_{n}^{\vee}\right\}
        \end{array}
\]
\end{definition}

\begin{uremark} 
  H\"aufig wird statt $\Gtlev<\Sp(\Lambda,\ZZ)$ mit der zu ihr
  konjugierten Gruppe $\Gamma_t=R_t^{-1}\Gtlev R_t<\Sp(4,\QQ)$,
  $R_t=\diag(1,1,1,t)$ gearbeitet. Diese hat den Vorteil, da\ss{} sie
  durch gebrochen-lineare Transformationen auf $\HH_2$ operiert. Bei
  der Berechnung des Titsgeb\"audes ist dieser Unterschied aber ohne
  Bedeutung. Allerdings werden unter Verwendung der Gruppe $\Gamma_t$
  statt $\Lambda$-isotroper Unterr\"aume $J$-isotrope Unterr\"aume
  klassifiziert. Wir werden im weiteren Verlauf isotrop als
  $\Lambda$-isotrop verstehen.
\end{uremark}

Zuerst werden wir noch einige zahlentheoretische \"Uberlegungen anstellen.
F\"ur eine positive ganze Zahl $n$ setzen wir
\[
        \nu(1) := 1,\ \nu(2) := 3 ~\mbox{und}~ \nu(n) := \frac{1}{2} n^{2}
\prod_{q|n,q \prim} (1-q^{-2})~\mbox{f\"ur}~n \ge 3
\]
und definieren die Menge der Nichttorsionselemente
\[
	\begin{array}{rcl}
\cN(n) & := & \left\{(a,b) \in \ZZ_{n}^{2} \mid \lambda (a,b) \neq 0 \mbox{ wenn } 0
  \neq \lambda \in \ZZ_{n} \right\}\\
	& =& \left\{(a,b)\in\ZZ_n^2\mid \ggT(n,a,b)=1\right\}
	\end{array}
\]
sowie
\[
\cM(n) := \cN(n)/\pm 1
\]
F\"ur jeden Teiler $k$ von $n$ sei 
\[
\cN_k(n)=\{(a,b)\in\ZZ_n^2\mid \ggT(n,a,b)=k\}.
\] 
Insbesondere ist $\cN_1(n)=\cN(n)$ und 
$\cN_n(n)=\{(0,0)\}$. Offensichtlich ist
$\ZZ_n^2$ die disjunkte Vereinigung $\ZZ_n^2=\bigcup\limits_{k|n}\cN_k(n)$, so
da{\ss} $n^2=\sum\limits_{k|n}\#\cN_k(n)$ ist. 
\begin{lemma}\label{AnzRest}
  Es sei $n \ge 2$. Die Anzahl der Restklassen in $\cN(n)$
  betr\"agt
\[
        n^2\prod_{q|n,q \prim}(1-q^{-2})
\]
\end{lemma}
{\it Beweis:\/}
\cite[Hilfssatz 1.2.3]{Zi}

Es sei $\ZZ_n^\times$ die Gruppe der Einheiten aus $\ZZ_n$. Dann operiert
$\ZZ_n^\times$ frei auf $\cN(n)$ durch $\lambda \cdot (v,w) :=
(\lambda v,\lambda w)$. Wir betrachten die Menge
\[
        \cO(n) := \cN(n)/\ZZ_n^\times 
\]
und setzen
\[
        \tilde{\nu}(1) = 1,\tilde{\nu}(2)= 3 ~\mbox{und}~ \tilde{\nu}(n) =
\frac{1}{\phi(n)}n^{2}\prod_{q|n,q \prim} (1-q^{-2})~\mbox{f\"ur}~n \ge 3.
\]
wobei $\phi$ die Eulersche $\phi$-Funktion ist. Aus Lemma \ref{AnzRest}
folgt unmittelbar, da\ss{} $\tilde{\nu}(n)$ die Anzahl der Restklassen in
$\cO(n)$ ($n\ge 2$) bestimmt. Mit $\tilde{\phi}(n)$ werden wir die Anzahl der
Elemente in $\ZZ_n^\times/\pm 1$ bezeichnen, also
\[
        \tilde{\phi}(1) = \tilde{\phi}(2) = 1 ~\mbox{und}~\tilde{\phi}(n) =
\frac{1}{2} \phi(t)~\mbox{f\"ur}~n \ge 3
\]
\begin{lemma}\label{laggT}
Es seien $x_1,x_2,x_3,x_4 \in \ZZ$ mit $\ggT(x_1,x_2,x_3,x_4)=1$ und
$x_2 \neq 0$. Dann gibt es $\lambda,\mu \in \ZZ$, so da\ss{}
$\ggT(x_1 + \lambda x_3 + \mu x_4,x_2) = 1$ ist.
\end{lemma}
{\it Beweis:\/}
\cite[Lemma~3.35]{HKW}

\section{Die Geb\"aude ${\cal T}(\Gt)$ und ${\cal T}(\Gtlev),~t$
quadratfrei}
%\markright{${\cal T}(\Gt)$ und ${\cal T}(\Gtlev)$} 

\begin{definition}
Es sei $v = (v_1,v_2,v_3,v_4)$ ein primitiver Vektor aus $\ZZ^4$. Dann
hei\ss{}t
\[
        d_{t}(v) := \ggT(v_1,v_3,t)
\]
der {\it $t$-Divisor} von $v$.
\end{definition}

Zur Berechnung der Titsgeb\"aude werden uns folgende Matrizen aus $\Gt$
(f\"ur $k\in\ZZ$) n\"utzlich sein, die durch Rechtsmultiplikation auf $\ZZ^4$ operieren:
\[
\begin{array}{cc}
 M_1(k) = \pmat{ 1 & k & 0 & 0 \cr
   0 & 1 & 0 & 0 \cr
   0 & 0 & 1 & 0 \cr
   0 & 0 & -tk & 1 } 
& M_2(k) = \pmat{ 1 & 0 & 0 & k \cr
   0 & 1 & tk & 0 \cr
   0 & 0 & 1 & 0 \cr
   0 & 0 & 0 & 1 } \\
 M_3(k) = \pmat{ 1 & 0 & 0 & 0 \cr
   -tk & 1 & 0 & 0 \cr
   0 & 0 & 1 & k \cr
   0 & 0 & 0 & 1 }
&
\end{array}
\]
und f\"ur $\pmat{a & b \cr c & d} \in \SL(2,\ZZ)$
\[
\begin{array}{lc}
 j_{1}  \left(\pmat {a & b \cr c & d}\right) =  
     \pmat{ a & 0 & b & 0 \cr
       0 & 1 & 0 & 0 \cr
       c & 0 & d & 0 \cr
       0 & 0 & 0 & 1 }  &
j_{2}\left(\pmat {a & b \cr c & d}\right) = 
     \pmat{ 1 & 0 & 0 & 0 \cr
       0 & a & 0 & b \cr
       0 & 0 & 1 & 0 \cr
       0 & c & 0 & d }
\end{array}
\]

Es gilt: $M_i(n)\in\Gtn \subseteq \Gtlev$,
$j_i\big(\Gamma_1(n)\big)\subset\Gtn$ und
$j_i\big(\Gamma_1(t)\big)\subset\Gtlev$ wobei $\Gamma_1(k):=\{\gamma
\in \SL(2,\ZZ)\mid \gamma \equiv \bbEins_2\mod k\}$.

Einen Vektor der Form $(0,a,0,b)$ werden wir im folgenden mit $v_{(a,b)}$ und
den Vektor $(1,0,0,0)$ mit $v_0$ notieren.

\begin{theorem}\label{GeradenBahnen}
  Zwei primitive Vektoren aus $\ZZ^4$, $v = (v_1,v_2,v_3,v_4)$ und $w
  =(w_1,w_2,w_3,w_4)$ sind genau dann $\Gt$-\"aquivalent, falls sie
den gleichen $t$-Divisor $r=d_t(v)=d_t(w)$ haben. Sie sind genau dann
  $\Gtlev$-\"aquivalent, falls auch $(v_2,v_4) \equiv (w_2,w_4) \mod r$ gilt.
\end{theorem}

{\it Beweis:\/} Der Beweis wird in 3 Schritten gef\"uhrt. Zun\"achst werden
wir einen beliebigen primitiven Vektor $v = (v_1,v_2,v_3,v_4) \in \ZZ^4$
mit $d_t(v) = r$ verm\"oge $\Gtlev$ in einen Vektor der Form
$(r,a_2,0,a_4)$ mit $\ggT(a_2,a_4)= 1$ transformieren, den wir mit $\Gt$ in
den Vektor $(r,1,0,0)$ \"uberf\"uhren werden. Im 2.Schritt wird die
Invarianz des $t$-Divisors eines primitiven Vektors unter $\Gt$ gezeigt und
letztlich wird bewiesen, da\ss{} zwei Vektoren $(r,v_2,0,v_4)$ und
$(r,w_2,0,w_4)$ genau dann $\Gtlev$-\"aquivalent sind, falls $(v_2,v_4)
\equiv (w_2,w_4) \mod r$ ist.

{\bf 1.Schritt: } Es sei also $(v_1,v_3) \equiv (0,0) \mod r$. Ohne
Einschr\"ankung sei $v_1 \neq 0$. Setze $x_1 = \frac{1}{r}v_3$, $x_2 =
\frac{1}{r}v_1 \neq 0$, $x_3 = \frac{t}{r}v_2$ und $x_4 = -\frac{t}{r}v_4$.
Da $r$ der $t$-Divisor von $v$ ist, gilt $\ggT(x_1,x_2,\frac{t}{r}) = 1$
und somit ist $\ggT(x_1,x_2,x_3,x_4) =
\ggT(x_1,x_2,\frac{t}{r}v_2,\frac{t}{r}v_4) = \ggT(x_1,x_2,v_2,v_4) = 1$,
da $v$ primitiv ist. Wenden wir Lemma \ref{laggT} auf $(x_1,x_2,x_3,x_4)$
an, so liefert dieses $\lambda,\mu \in \ZZ$ mit $\ggT(x_2,x_1 + \lambda x_3
+ \mu x_4) = 1$. Unter $M_1(\mu)M_2(\lambda)$ werden die Eintr\"age
$(v_1,v_3)$ nach $(v_1,v_3 +\lambda t v_2 - \mu t v_4 + \lambda \mu t
v_1)$ transformiert und es gilt
\begin{eqnarray*}
\ggT(v_1,v_3 + \lambda t v_2 - \mu t v_4 + \lambda \mu t v_1) &=& \ggT(v_1,v_3
+ \lambda t v_2 - \mu t v_4)\\
 &=& \ggT(r x_2,r x_1 + \lambda r x_3 + \mu r
x_4)\\
 &=& r \ggT(x_2,x_1 + \lambda x_3 + \mu x_4)\\
 &=& r
\end{eqnarray*}
Wir k\"onnen also $\ggT(v_1,v_3) = r$ annehmen. Ein geeignetes Element aus
$j_1(\SL(2,\ZZ))$ bringt diesen nach $(r,v_2,0,v_4)$. Ist nun $v_2 = 0$, so
transformiert $M_1(1)$ den Vektor $v$ nach $(r,r,-tv_4,v_4)$, der wiederum
durch $j_1(\SL(2,\ZZ))$ auf $(r,r,0,v_4)$ geht. Wir k\"onnen also auch $v_2
\neq 0$ annehmen. Nun ist $\ggT(r,v_2,v_4) = 1$, da die Eigenschaft
primitiv eines Vektors aus $\ZZ^4$ unter $\Gt$ erhalten bleibt. Wenden wir
nun Lemma \ref{laggT} auf $(v_4,v_2,r,0)$ an, so liefert uns diese ein $\mu
\in \ZZ$ mit $\ggT(v_2,v_4 + \mu r) = 1$. Anwendung der Matrix $M_2(\mu)$
transformiert $v$ nach $(r,v_2,\mu t v_2,v_4 + \mu r)$, der via
$j_1(\SL(2,\ZZ))$ auf $(r,v_2,0,v_4+ \mu r)$ geht. Ein beliebiger
primitiver Vektor $v$ mit $d_t(v) = r$ kann also nur unter Verwendung von
Matrizen aus $\Gtlev$ auf einen Vektor der Form $(r,a_2,0,a_4)$ gebracht
werden. Da $\ggT(v_2,v_4 + \mu r) = 1$ ist, finden wir letztlich ein Element
aus $j_2(\SL(2,\ZZ))$, das $(r,v_2,0,v_4 + \mu r)$ auf $(r,1,0,0)$
transformiert.

{\bf 2.Schritt: } Nach dem 1.Schritt reicht es aus zu zeigen, da\ss{} ein
Vektor der Form $v = (r,1,0,0) \in \ZZ^4, r|t$ die Eigenschaft $d_t(v) = r$
unter der Operation von $\Gt$ beh\"alt. Als Verallgemeinerung der
\"Uberlegungen in \cite[Lemma~0.5]{HW} ergeben sich f\"ur eine Matrix $(m_{ij})
\in \Gt$ notwendigerweise die Kongruenzen
\begin{eqnarray*}
        m_{21} \equiv m_{41} \equiv m_{23} \equiv m_{43} &\equiv & 0 \mod t \\
        m_{11}m_{33} - m_{13}m_{31} &\equiv & 1 \mod t
\end{eqnarray*}
Wenden wir also eine Matrix $M = (m_{ij}) \in \Gt$ auf $v$ an, so erhalten wir
\[
 v \cdot M = (rm_{11} + m_{21},\ast,rm_{13} + m_{23},\ast)
\] 
Der $t$-Divisor von $v \cdot M$ ist
\begin{eqnarray*}
d_t(v\cdot M)&=&\ggT(rm_{11} + m_{21},rm_{13} + m_{23},t)\\
             &=& \ggT(rm_{11},rm_{13},t)\\
             &=& r \ggT(m_{11},m_{13},\frac{t}{r}).
\end{eqnarray*}
W\"are $\ggT(m_{11},m_{13},\frac{t}{r}) \neq 1$, dann auch
$\ggT(m_{11},m_{13},t) \neq 1$ im Widerspruch zu $m_{11}m_{33} -
m_{13}m_{31} \equiv 1 \mod t$.

{\bf 3.Schritt: } Es ist klar, da\ss{} zwei $\Gtlev$-\"aquivalente
primitive Vektoren $v = (v_1,v_2,v_3,v_4)$ und $w = (w_1,w_2,w_3,w_4)$ die
Kongruenz $(v_2,v_4) \equiv (w_2,w_4) \mod r$
erf\"ullen m\"ussen. Hier geht im wesentlichen ein, da\ss{} f\"ur eine
Matrix $M = (m_{ij}) \in \Gtlev$ notwendigerweise die Kongruenz
\[
m_{22}-1 \equiv m_{44} -1 \equiv m_{24} \equiv m_{42} \equiv 0 \mod
t~~(\dagger)
\]
gelten mu\ss{}, die nat\"urlich auch modulo $r$ gelten.
Ist nun $v \cdot M =w$ f\"ur ein $M = (m_{ij}) \in \Gtlev$, dann ist
\begin{eqnarray*}
        w_2 &=& m_{12}v_1 + m_{22}v_2 + m_{32}v_3 + m_{42}v_4 \\
        w_4 &=& m_{14}v_1 + m_{24}v_2 + m_{34}v_3 + m_{44}v_4
\end{eqnarray*}
also 
\begin{eqnarray*}
        w_2 \equiv m_{12}v_1 + v_2 + m_{32}v_3 \mod r \\
        w_4 \equiv m_{14}v_1 + m_{34}v_3 + v_4 \mod r
\end{eqnarray*}
wegen $(\dagger)$ und damit $(v_2,v_4) \equiv (w_2,w_4) \mod r$, weil $v_1
\equiv v_3 \equiv 0 \mod r$. Die folgenden \"Uberlegungen zeigen, da\ss{}
diese Kongruenzbeziehungen auch hinreichend sind. Nach dem 1.Schritt
k\"onnen wir ohne Einschr\"ankung annehmen, da\ss{} die Vektoren $v$ und
$w$ durch $(r,v_2,0,v_4)$ bzw. $(r,w_2,0,w_4)$ gegeben sind. Es ist
$(v_2,v_4) \equiv (w_2,w_4) \mod r$, also $v_2 = kr + w_2$ und $v_4 = lr +
w_4$ f\"ur geeignete $k,l \in \ZZ$. Der Vektor $w$ geht unter Anwendung der
Matrix $M_1(k)M_2(l)$ auf $(r,w_2+kr,ltw_2 - ktw_4-kltr, w_4 + lr)$, der
durch ein geeignetes Element aus $j_1(\SL(2,\ZZ)) \subset \Gtlev$ auf $v$
geht, da $\ggT(r,ltw_2-ktw_4-kltr) = r$ ist.  \qed

\begin{cor}\label{gtlevtrans}
Die Gruppe $\Gtlev$ operiert transitiv auf der Menge der primitiven Vektoren
mit $t$-Divisor~$1$.
\end{cor}

Damit k\"onnen wir nun die Geraden in $\QQ^{4}$ klassifizieren. Bezeichnen
wir mit $\mu(t)$ die Anzahl der Teiler von $t$, so gilt
\begin{theorem}
  Die eindimensionalen (und somit isotropen) Unterr\"aume in $\QQ^4$
  zerfallen unter der Operation von $\Gt$ in genau $\mu(t)$ Bahnen.
  
  Es seien $r_i, i=1,\ldots, \mu(t)$ die Teiler von $t$. Dann ist die Anzahl
  der $\Gtlev$-Bahnen eindimensionaler Unterr\"aume durch
\[
        \psi (t) := \sum_{i=1}^{\mu(t)}\nu(r_i) = \left\{ \begin{array}{c}
(t^2+1)/2 \cr (t^2+4)/2 \end{array}~ \mbox{falls}~ \begin{array}{l} 
t~\mbox{ungerade} \cr t~\mbox{gerade}\end{array} \right.
\]
bestimmt.
\end{theorem}
{\it Beweis:\/} Zuerst wird jede Gerade in $\QQ^{4}$ mit dem bis auf
Vorzeichen eindeutig bestimmten primitiven Vektor aus $\ZZ^{4}$, der die
Gerade erzeugt, identifiziert. Nach Satz \ref{GeradenBahnen} sind zwei
primitive Vektoren genau dann $\Gt$-\"aquivalent, wenn sie den gleichen
$t$-Divisor haben. F\"ur $t$ ergeben sich hieraus $\mu(t)$ M\"oglichkeiten.

Unter der Operation von $\Gtlev$ bleibt zus\"atzlich die Auswahl von
$(a_2,a_4)\mod r$ und modulo Vorzeichen, das hei\ss{}t die Auswahl eines
Element aus $\cM(r)$, deren M\"achtigkeit wir in Lemma \ref{AnzRest} mit
$\nu(r)$ bestimmt haben.\\
Es ist leicht zu sehen, da{\ss} die Abbildung $\cN_r(t)\to \cN_1(t/r)$ gegeben
durch $(a,b)\mapsto (a/r, b/r)$ bijektiv ist. Daraus folgt sofort
\[
\begin{array}{rcl}
t^2&=&\sum\limits_{r|t}\#\cN(r)\\
   &=& \sum\limits_{r|t, r<3}\#\cN(r)+\sum\limits_{r|t, r\geq
3}\#\cN(r)\\
   &=& \sum\limits_{r|t, r<3}\#\cN(r)+\sum\limits_{r|t, r\geq
3}2\nu(r)\\
&=& -\sum\limits_{r|t, r<3}\#\cN(r)+\sum\limits_{r|t}2\nu(r)\\
&=& -\sum\limits_{r|t, r<3}\#\cN(r)+2\psi(t)
\end{array}
\]
Der erste Term ist $-1$ f\"ur $t$ ungerade und $-1-3=4$ f\"ur $t$ gerade.
\qed

Kommen wir nun zu der Klassifizierung der isotropen Ebenen. Erst hier wird
die Forderung quadratfrei an $t$ zum Tragen kommen. Jede isotrope Ebene $h$
aus $\QQ^4$ k\"onnen wir als Erzeugnis zweier primitiver Vektoren $v$ und
$w$ schreiben. Ohne Einschr\"ankung werden wir zus\"atzlich an die
erzeugenden Vektoren die Forderung stellen, da\ss{} sie das Gitter $h_{\ZZ}
:= h \cap \ZZ^4$ erzeugen.
\begin{lemma}\label{lateilerfremd}
  Es sei $t$ quadratfrei und $h = v \wedge w$ eine Ebene mit $d_t(v) = r,
  d_t(w) =s$ und $h_{\ZZ} = \ZZ v \oplus \ZZ w$. Dann sind $r$ und $s$
  teilerfremd.
\end{lemma}

{\it Beweis:\/}
Nach den obigen \"Uberlegungen k\"onnen wir ohne Einschr\"ankung $v =
(r,1,0,0)$ voraussetzen. Angenommen, $r$ und $s$ sind nicht teilerfremd.

Es sei $m:=\ggT(r,s)$. Wir schreiben $mr' = r$ und $ms' = s$ mit $\ggT(r',s')
=1$. Aus der Isotropieeigenschaft von $h$ ergibt sich $w_3 =
-\frac{t}{r}w_4$. F\"ur den $t$-Divisor von $w$ gilt $d_t(w) = \ggT(w_1,w_3,t)
= s$ und insbesondere, da $m$ die Zahl $s$ teilt
\[
        -\frac{t}{r}w_4 = w_3 \equiv w_1 \equiv 0 \mod m
\]
Nun sind aber $\frac{t}{r}$ und $m$ teilerfremd: H\"atten $\frac{t}{r}$ und
$m$ einen gemeinsamen Teiler $l>1$, so k\"onnen wir $\frac{t}{r} = t'l$ und
$m= m'l$ mit $\ggT(t',m') = 1$ schreiben. Dann w\"are aber $t = t'lr =
t'lmr' = tl^2m'r'$ im Widerspruch zu $t$ quadratfrei. Die Zahl $m$ teilt
also $r,w_1,w_3$ und $w_4$. Da die Vektoren $v$ und $w$ das Gitter
$h_{\ZZ}$ erzeugen, wird dieses auch von $w' = w - w_2 v =
(w_1-rw_2,0,w_3,w_4)$ und $v$ erzeugt. Insbesondere mu\ss{} $w'$ auch
primitiv sein, aber $m$ ist gemeinsamer Teiler von $w_1-rw_2,w_3$ und
$w_4$. \qed
\begin{theorem}\label{ladt1}
  Es sei $t$ quadratfrei und $h = v \wedge w$ eine isotrope Ebene mit
  $d_t(v) = r, d_t(w) = s$ und $h_{\ZZ} = \ZZ v \oplus \ZZ w$. Dann gibt es
  primitive Vektoren $\hat{v},\hat{w}$, so da\ss{} $h_{\ZZ} = \ZZ \hat{v}
  \oplus \ZZ \hat{w}$ und $d_t(\hat{v}) = 1$.
\end{theorem}
{\it Beweis:\/} Es sei $h = v \wedge w$ mit obigen Voraussetzungen gegeben.
Dann ist nach Lemma \ref{lateilerfremd} $\ggT(r,s) = 1$. F\"ur den Fall,
da\ss{} einer der Vektoren den $t$-Divisor $t$ hat, ist also nichts zu
zeigen.
 
Es sei $m = \min\{d_t(u) : u \in h_{\ZZ}\}$ und $\tilde{v}$ ein primitiver
Vektor mit $d_t(\tilde{v}) = m$. Da $v$ und $w$ das Gitter $h_{\ZZ}$
aufspannen, gibt es ganze Zahlen $\lambda,\mu \in \ZZ$, so da\ss{}
$\tilde{v} = \lambda v + \mu w$ gilt. Da $\tilde{v} $ primitiv ist, sind
insbesondere $\lambda$ und $\mu$ teilerfremd. W\"ahle $\nu,\xi \in \ZZ$ so,
da\ss{} $\lambda \nu + \mu \xi = 1$ ist und definiere $\tilde{w} := -\xi v
+ \nu w$. Die Vektoren $\tilde{v}$ und $\tilde{w}$
spannen damit das Gitter $h_{\ZZ}$ auf.

Angenommen, es gilt $m > 1$. Nach Satz \ref{GeradenBahnen} gibt es ein
$\gamma \in \Gt$, so da\ss{} $\tilde{v} \cdot \gamma = (m,1,0,0) =:
v'$ ist.  Setze $(w_1',w_2',w_3',w_4') = w' := \tilde{w} \cdot \gamma$
und $h' = v' \wedge w'$. Das Gitter $h'_{\ZZ}$ wird von $v'$ und $w'$
aufgespannt und es gilt $m = \min\{d_t(u')\mid u' \in
h_{\ZZ}'\}$. G\"abe es n\"amlich einen primitiven Vektor $u'$ mit
$d_t(u') < m$, dann w\"are $u'\cdot \gamma^{-1}$ Teil einer Basis von
$h_{\ZZ}$ und $d_t(u'\cdot \gamma^{-1}) < m$.

Wir setzen $l=\ggT(m,w'_1)$ und werden folgende F\"alle unterscheiden:

i) $l=m$. Ohne Einschr\"ankung sei $w_1' = 0$ (Ist n\"amlich $ w_1' \neq
0$, dann gibt es eine Zahl $k \in \ZZ$, so da\ss{} $w_1' = km$. Das Gitter
$h_{\ZZ}'$ wird dann auch von $v'$ und $w'-kv'$ erzeugt). Das Gitter
$h_{\ZZ}'$ wird von $v'$ und $w'$, also auch von $\tilde{v} := v' + w' =
(m,w_2' + 1,w_3',w_4')$ und $w'$ erzeugt. Insbesondere ist $\tilde{v}$
primitiv. Der $t$-Divisor von $w'$ ist $d_t(w') = \ggT(0,w_3',t) =: l'$. Da
$\tilde{v}$ und $w'$ das Gitter $h_{\ZZ}'$ aufspannen, gilt nach Lemma
\ref{lateilerfremd}, da\ss{} $\ggT(m,l')= 1$ ist. Dann ist $d_t(\tilde{v})
= \ggT(m,w_3',t) = \ggT(m,l') = 1$ im Widerspruch zu $m > 1$.

ii) $l<m$. Schreibe $w_1' = kl$ und $m = m'l$. Da $l = \ggT(m,w_1')$, gibt
es $\lambda',\mu' \in \ZZ$, so da\ss{} $\lambda'm + \mu'w_1' = l$, das
hei\ss{}t $\lambda'm' + \mu'k = 1$ und somit spannen die Vektoren
$\tilde{v} := \lambda' v' + \mu' w' = (l,\lambda' +
\mu'w_2,\mu'w_3',\mu'w_4')$ und $\tilde{w} := -kv' + m'w'$ das Gitter
$h_{\ZZ}'$ auf. Es gilt jedoch $d_t(\tilde{v}) = \ggT(l,\mu w_3',t) \le l <
m$ im Widerspruch zur Minimaleigenschaft von $m$. \qed
\begin{cor}\label{corerzeugnis}
  Jede isotrope Ebene $h$ kann als Erzeugnis zweier Vektoren $v$ und $w$
  mit $h_{\ZZ} = \ZZ v \oplus \ZZ w$, $d_t(v) = 1$ und $d_t(w) = t$
  geschrieben werden.
\end{cor}
{\it Beweis:\/} Nach Lemma \ref{ladt1} k\"onnen wir ohne
Einschr\"ankung annehmen, da\ss{} $d_t(v) = 1$ gilt. Dann gibt es nach
Satz \ref{GeradenBahnen} ein $\gamma$ aus $\Gt$, so da\ss{}
$v\cdot \gamma = v_0$ ist. Setze $\tilde{w} := w \cdot \gamma$ und $\tilde{h}
:= v_0 \wedge \tilde{w}$. Das Gitter $\tilde{h}_{\ZZ}$ wird von $v_0$ und
$\tilde{w}$ erzeugt. F\"ur $\tilde{w} =
(\tilde{w}_1,\tilde{w}_2,\tilde{w}_3,\tilde{w}_4)$ folgt aus der
Isotropieeigenschaft, da\ss{} $\tilde{w}_3 = 0$ ist. Das Gitter
$\tilde{h}_{\ZZ}$ wird auch von $w' = \tilde{w} - \tilde{w}_1v_0$ und
$v_0$ erzeugt. Der $t$-Divisor von $w'$ ist $d_t(w') = \ggT(0,0,t) =
t$. Nun ist $h = (v_0 \cdot \gamma^{-1}) \wedge (w'\cdot \gamma^{-1})$ und $h_{\ZZ} =
\ZZ (v_0 \cdot\gamma^{-1}) \oplus \ZZ (w'\cdot\gamma^{-1})$. Aufgrund der
Invarianz des $t$-Divisors unter $\Gt$ folgt die Behauptung.  \qed

\begin{theorem}
  Die Gruppe $\Gt$ operiert transitiv auf der Menge der iso\-tropen Ebenen.
  Die Menge der $\Gtlev$-Bahnen isotroper Ebenen wird von $\cO(t)$
  indiziert.
\end{theorem}
{\it Beweis:\/} Es sei $h = v \wedge w$ mit $h_{\ZZ} = \ZZ v \oplus \ZZ w$
und $d_t(v) = 1$ eine isotrope Ebene. Nach Satz \ref{GeradenBahnen} gibt
es ein $\gamma \in \Gt$, so da\ss{} $v \cdot \gamma = v_0$ ist. Wir setzen
$\tilde{w} := w \cdot \gamma$ und $\tilde{h} = v_0 \wedge \tilde{w}$. Ist
$\tilde{w} = ({\tilde w}_1,{\tilde w}_2, {\tilde w}_3,{\tilde w}_4)$, so
folgt aus der Isotropieeigenschaft $\tilde{w}_3 =0$. Indem wir
$\tilde{w}_1$ durch $\tilde{w} - \tilde{w}_1v_0$ ersetzen, k\"onnen wir
$\tilde{w}_1 = 0$ voraussetzen. Da $v_0$ und $\tilde{w}$ das Gitter
$\tilde{h}_{\ZZ}$ erzeugen, ist $\ggT(\tilde{w}_2,\tilde{w}_4) = 1$.
W\"ahlen wir nun $\lambda,\mu \in \ZZ$ mit $\lambda \tilde{w}_2 + \mu
\tilde{w}_4 = 1$, so ist $\gamma =\pmat{\lambda & -\tilde{w}_4 \cr \mu &
\tilde{w}_2} \in \SL(2,\ZZ)$ und $j_{2}(\gamma) \in \Gt $ \"uberf\"uhrt die
Ebene $\tilde{h}$ nach $v_0 \wedge v_{(1,0)}$, womit die erste Aussage
gezeigt ist.

Es sei nun $h = v \wedge w$ wie oben. Da $\Gtlev$ auf der Menge der
primitiven Vektoren mit $t$-Divisor 1 transitiv operiert, k\"onnen wir ohne
Einschr\"ankung $h = v_0 \wedge w$ annehmen. Wieder k\"onnen wir annehmen,
da\ss{} $w_1 = 0$ und aufgrund der Isotropieeigenschaft von $h$ ist $w_3 =
0$. Die Eintr\"age $w_2$ und $w_4$ m\"ussen notwendigerweise teilerfremd
sein, da $w$ sonst nicht primitiv w\"are. Jede isotrope Ebene ist also
\"aquivalent zu einer Ebene der Form $h = v_0 \wedge v_{(w_2,w_4)}$ mit
$\ggT(w_2,w_4)=1$. Das Paar $(w_2,w_4)$ kann kein Torsionselement in
$\ZZ_{t}^2$ sein und es definiert deshalb eine mit
$\left[w_2,w_4\right]_{\cO(t)}$ bezeichnete Klasse in~$\cO(t)$. Wir werden
zeigen, da\ss{} die Zuordnung
\[
        \Phi : \left[{h}\right]_{\Gtlev} \mapsto
\left[w_2,w_4\right]_{\cO(t)}
\]
wohldefiniert und bijektiv ist.

i) $\Phi$ is wohldefiniert:\\
Dazu sei $\tilde{h} := v_0 \wedge v_{(\tilde{w}_2,\tilde{w}_4)}$
eine zu $h$ \"aquivalente isotrope Ebene, also $h\cdot \tilde{\gamma} =
\tilde{h}$ f\"ur ein $\tilde{\gamma} \in \Gtlev$. Dann wird das Gitter
$\tilde{h}_{\ZZ}$ von $v_0 \cdot \tilde{\gamma}$
und $w\cdot\tilde{\gamma}$ aufgespannt, das hei\ss{}t
\[
        \begin{array}{rcl}v_0 \cdot \tilde{\gamma} &  = & av_0 +
b v_{(\tilde{w}_2,\tilde{w}_4)} \\
        w \cdot \tilde{\gamma} & = & cv_0 + dv_{(\tilde{w}_2,\tilde{w}_4)}
        \end{array}
\]
f\"ur ein $\pmat{a & b \cr c & d} \in \GL(2,\ZZ)$. Der $t$-Divisor ist
unter der Operation von $\Gtlev$ invariant, es gilt damit
\[
        d_t(w\cdot \tilde{\gamma}) = d_t(w) = \ggT(0,0,t) = t
\]
und somit mu\ss{} $c \equiv 0\mod t$ gelten. Da $c$ und $d$ teilerfremd
sind, gilt $\ggT(d,t) = 1$. Die Zahl $d$ repr\"asentiert also eine Einheit
in $\ZZ_t$. Schreiben wir $\tilde{\gamma}(w) =
(\hat{w}_1,\hat{w}_2,\hat{w}_3,\hat{w}_4)$, so ist also
$(\hat{w}_2,\hat{w}_4) \equiv (d\tilde{w}_2,d\tilde{w}_4) \mod t$. Nun ist
die Restklasse $(w_2,w_4)\mod t$ unter $\Gtlev$ invariant, es gilt also
\[
(w_2,w_4) \equiv (\hat{w}_2,\hat{w}_4)\equiv
(d\tilde{w}_2,d\tilde{w}_4)\mod t
\]
Die Restklassen $(w_2,w_4) \mod t$ und $(\tilde{w}_2,\tilde{w}_4)\mod t$
stimmen bis auf ein Vielfaches eines Repr\"asentanten einer Einheit in
$\ZZ_t$ \"uberein und somit ist $\left[w_2,w_4\right]_{\cO(t)} =
\left[\tilde{w}_2,\tilde{w}_4\right]_{\cO(t)}$.

ii) $\Phi$ ist surjektiv:\\
F\"ur eine Restklasse $\left[w_2,w_4\right]_{\cO(t)}$ w\"ahlen wir einen
Repr\"asentanten $(w'_2,w'_4)$ mit $\ggT(w'_2,w'_4) = 1$. Dann ist die Bahn
von isotropen Ebenen, die durch $h = v_0 \wedge v_{(w'_2,w'_4)}$
repr\"asentiert wird, im Urbild von $\left[w_2,w_4\right]_{\cO(t)}$.

iii) $\Phi$ ist injektiv:\\
Es seien $(w_2,w_4)$ und $(w'_2,w'_4)$ zwei Repr\"asentanten aus
$\left[w_2,w_4\right]_{\cO(t)}$. Dann ist $(ew_2,ew_4) \equiv
(w'_2,w'_4)\mod t$ f\"ur einen Repr\"asentanten $e$ einer Einheit in
$\ZZ_t$, also insbesondere $\ggT(e,t) = 1$. W\"ahle nun $\lambda,\mu \in
\ZZ$ so, da\ss{} $\pmat{\mu & -\lambda \cr t & e}$ aus $\SL(2,\ZZ)$ ist.
Die Vektoren
\[
        \begin{array}{rcl}      
\tilde{v} & := & \mu v_0 - \lambda v_{(w_2,w_4)} \cr
        \tilde{w} & := & tv_0 + ev_{(w_2,w_4)}
        \end{array}
\]
spannen dann das Gitter $h_{\ZZ}$ auf.

Es gilt $d_t(\tilde{v})= \ggT(\mu,0,t) = 1$ und $d_t(\tilde{w}) =
\ggT(t,0,t) = t$. Nach Korollar~\ref{gtlevtrans} gibt es ein $\gamma \in
\Gtlev$ mit $\tilde{v}\cdot\gamma = v_0$. Wir setzen $\hat{w} :=
\tilde{w} \cdot \gamma$ und schreiben $\hat{w} =
(\hat{w}_1,\hat{w}_2,\hat{w}_3,\hat{w}_4)$. Dann mu\ss{} wieder $\hat{w}_3
= 0$ sein und ohne Einschr\"ankung k\"onnen wir $\hat{w}_1 = 0$ annehmen.
Da die Kongruenzklasse $(ew_2,ew_4)\mod t$ invariant unter $\Gtlev$ ist,
gilt
\[
        (\hat{w}_2,\hat{w}_4) \equiv (ew_2,ew_4) \equiv (w_2',w_4')\mod t
\]
Nach \cite[Lemma~1.42]{Sh} gibt es ein $g \in \Gamma_1(t)$, so da\ss{}
$(\hat{w}_2,\hat{w}_4)\cdot g = (w_2',w_4')$ ist. Die Matrix $\tilde{\gamma}
:= j_2(g)$ ist aus $\Gtlev$ und
\[
        (v_0 \cdot \tilde{\gamma}) \wedge (v_{(\hat{w}_2,\hat{w}_4)}\cdot \tilde{\gamma})=v_0\wedge
v_{(w_2',w_4')} = h'
\]
Das hei\ss{}t also, $\gamma\tilde{\gamma} \in \Gtlev$ \"uberf\"uhrt
die isotrope Ebene $h$ nach $h'$.  \qed

Damit sind die Ecken der Titsgeb\"aude ${\cal T}(\Gtlev)$ und ${\cal
  T}(\Gt)$ vollst\"andig be\-schrieben. Was bleibt, ist die Untersuchung
der Kanten der Titsgeb\"aude.  Im Fall ${\cal T}(\Gt)$ k\"onnen wir
uns auf die Ebene $(1,0,0,0) \wedge (0,1,0,0)$ beschr\"anken, die bis
auf \"Aquivalenz die Unterr\"aume enth\"alt, die von primitiven
Vektoren $(s_i,1,0,0)$ aufgespannt werden, wobei $s_i$ alle Teiler von
$t$ durchl\"auft. F\"ur ${\cal T}(\Gt)$ ergibt sich also eine
$(\mu(t)_1,1_{\mu(t)})$-Konfiguration:

\begin{figure}[h]
  \begin{center}
\epsfbox{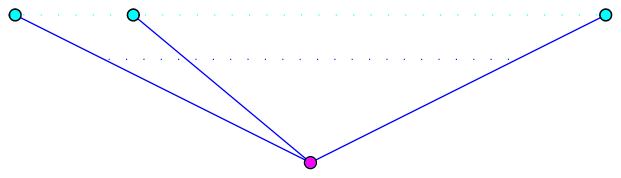} 
  \end{center}
\begin{center}
Bild 1: ${\cal T}(\Gt)$
\end{center}
\end{figure}

\begin{lemma}
  Sind $r_i,i=1,\ldots ,\mu(t)$ die Teiler von $t$, so enth\"alt jede
  isotrope Ebene bis auf $\Gtlev$-\"Aquivalenz genau
\[
        \sum_{i=1}^{\mu(t)}\tilde{\phi}(r_i)
\]
Geraden.
\end{lemma}
{\it Beweis:\/} Es reicht aus, dazu die isotrope Ebene $h = v_0 \wedge
v_{(1,0)}$ zu betrachten. Diese enth\"alt bis auf \"Aquivalenz die Geraden,
die von primitiven Vektoren der Form $(r_i,k,0,0)$ aufgespannt werden. Was
bleibt, ist die Auswahl von $k$ aus $\ZZ_{r_i}^{\times}/\pm 1$, dessen
M\"achtigkeit wir mit $\tilde{\phi}(r_i)$ bestimmt haben. \qed

Wir bezeichnen mit $\gh{w_2,w_4}$ die Bahn von isotropen Ebenen, die durch
$\left[w_2,w_4\right]_{\cO(t)}$ bzw. mit $\gl{x_2,x_4}$ die Bahn von
Geraden, die durch den Teiler $r$ und der Restklasse
$\left[x_2,x_4\right]_{\cM(r)}$ bestimmt ist.

\begin{lemma}
  Je zwei Ecken $\gh{w_2,w_4}$ und $\gl{x_2,x_4}$ im Titsgeb\"aude ${\cal
  T}(\Gtlev)$ werden genau dann miteinander verbunden, wenn $(w_2,w_4)$ und
  $(x_2,x_4)$ die selben Restklassen in $\cO(r)$ bestimmen.
\end{lemma}
{\it Beweis:\/} Da $\left[w_2,w_4\right] \in \cO(t)$ ist,
repr\"asentiert $(w_2,w_4)$ kein Torsionselement in $\ZZ^2_t$. Die
Zahl $r$ teilt $t$ und damit ist $(w_2,w_4)$ auch kein Repr\"asentant
eines Torsionselements in $\ZZ^2_r$, das hei\ss{}t,
$\left[w_2,w_4\right]_{\cO(r)}$ ist wohldefiniert.

($\Rightarrow$)\\
 Seien $\gh{w_2,w_4}$ und $\gl{x_2,x_4}$ miteinander
verbunden. Dann gibt es f\"ur jedes Element $h$ aus $\gh{w_2,w_4}$
einen Repr\"asentanten $\ell\in\gl{x_2,x_4}$ mit $\ell \subset h$. Wir
w\"ahlen $h = v_0 \wedge v_{(w_2,w_4)}$. Mit $v = (v_1,v_2,v_3,v_4)$
bezeichnen wir den bis auf Vorzeichen eindeutig bestimmten primitiven
Vektor, der den Unterraum $\ell$ aufspannt. Dann l\"a\ss{}t sich $v$
als Linearkombination $v = av_0 + bv_{(w_2,w_4)}$ f\"ur $a,b \in \ZZ$
schreiben. Da $\ell$ Repr\"asentant der Bahn $\gl{x_2,x_4}$ ist, gilt
$d_t(v) = r$ und $(v_2,v_4) \equiv (x_2,x_4) \mod r$. Also mu\ss{} $a$
die Kongruenz $a \equiv 0 \mod r$ erf\"ullen. Dann ist $b$ aber
Repr\"asentant einer Einheit in $\ZZ_r$, weil $v$ sonst nicht primitiv
w\"are. Somit stimmen $(v_2,v_4)$ und $(w_2,w_4)$ bis auf eine Einheit
in $\ZZ_r$ \"uberein und damit gilt $\left[w_2,w_4\right]_{\cO(r)} =
\left[v_2,v_4\right]_{\cO(r)} = \left[x_2,x_4\right]_{\cO(r)}$.

($\Leftarrow$)\\
 Es gelte $\left[x_2,x_4\right]_{\cO(r)} =
\left[w_2,w_4\right]_{\cO(r)}$ und $h = v_0 \wedge
v_{(w_2,w_4)}\in\gh{w_2,w_4}$.  Nun ist $\left[x_2,x_4\right]_{\cM(r)}
= \left[ew_2,ew_4\right]_{\cM(r)}$ f\"ur einen Repr\"asentanten $e$
einer Einheit in $\ZZ_r$.  Betrachte den Unterraum $\ell$, der von
$(r,ew_2,0,ew_4)=rv_0+ev_{(w_2,w_4)}\in h$ erzeugt wird. Dann gilt
$\ell\subset h$, und nach $\left[ew_2,ew_4\right]_{\cM(r)} =
\left[x_2,x_4\right]_{\cM(r)}$ auch $\ell\in\gl{x_2,x_4}$. \qed

Mit diesen Angaben k\"onnen wir letzlich das Titsgeb\"aude ${\cal T}(\Gtlev)$
angeben:
\begin{theorem}
  Das Titsgeb\"aude ${\cal T}(\Gtlev)$ ist ein Graph mit $\psi(t)$ Ecken,
  die von Geraden $\ell_0 = \QQ v_0$ und $\ell_{r_i,(v_2,v_4)} =
  \QQ(r_i,v_2,0,v_4)$ repr\"asentiert werden, wobei $r_i$ alle Teiler $\neq 1$ von
  $t$ und $(v_2,v_4)$ alle Restklassen in $\cM(r_i)$ durchl\"auft, sowie
  $\tilde{\nu}(t)$ Ecken, die von Ebenen $h_{\left[w_2,w_4\right]} = v_0
  \wedge v_{(w_2,w_4)}$ repr\"asentiert werden, wobei $(w_2,w_4)$ alle
  Restklassen in $\cO(t)$ durchl\"auft. Jede Ecke $h_{\left[w_2,w_4\right]}$
  wird mit der Ecke $\ell_0$ verbunden und $\ell_{r_i,(v_2,v_4)}$ wird genau dann mit der Ecke
  $h_{\left[w_2,w_4\right]}$ verbunden, wenn $(v_2,v_4)$ und $(w_2,w_4)$
  Repr\"asentanten der gleichen Restklasse in $\cO(r_i)$ sind.
\end{theorem}

\begin{figure}[htbp]
  \begin{center}
\epsfbox{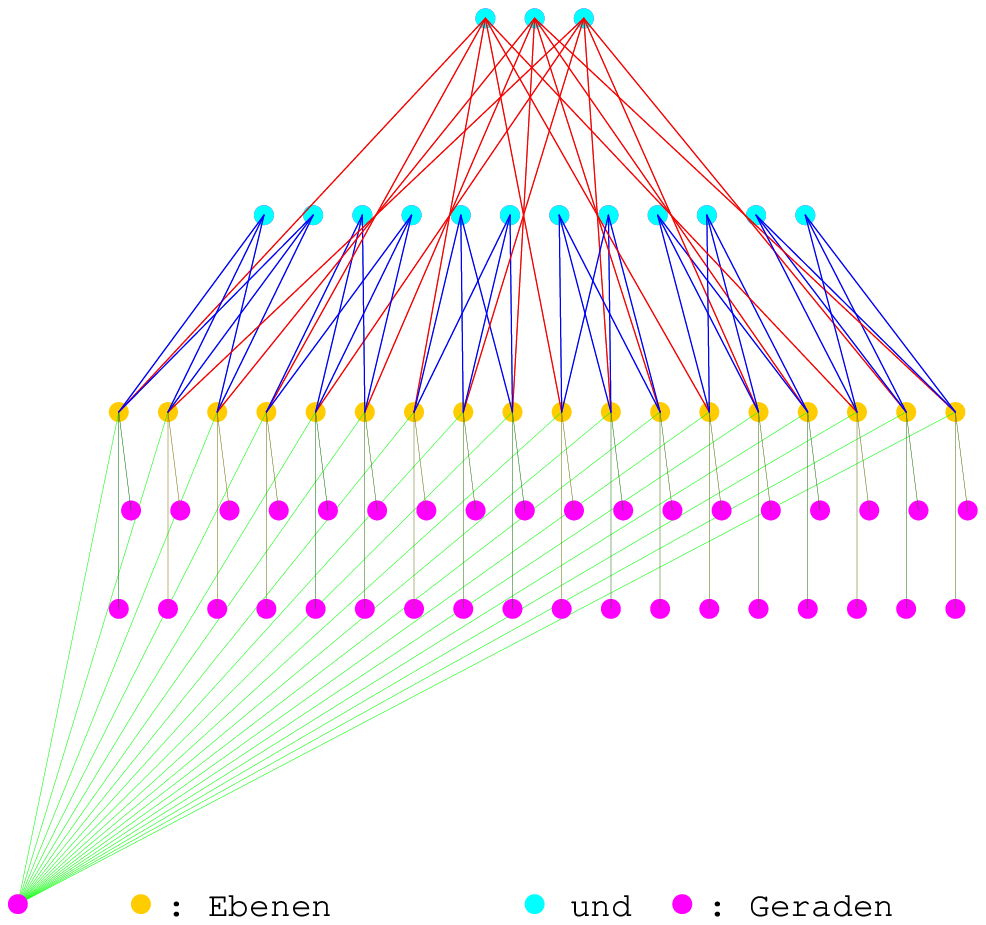} 
  \end{center}
\begin{center}
Bild 2: ${\cal T}({{\tilde{\Gamma}_{1,10}}})$
\end{center}
\end{figure}

\section{Die Geb\"aude ${\cal T}(\Gpq)$, $p,q$ prim, $p \geq 3$.}
%\markright{${\cal T}(\Gpq)$}
Wir werden das Titsgeb\"aude der Paramodulgruppe $\Gpq$ zur Stufe~$q$ und
Typ~$(1,p)$ berechnen. Der Einfachheit halber setzen wir voraus, da\ss{}
$p$ sowie $q$ Primzahlen sind, und $p\ge 3$ ist. Die entsprechenden
Mo\-dulr\"aume sind (noch) nicht so gr\"undlich untersucht worden: wir
erw\"ahnen aber den Fall $p=3$, $q=2$, das der Nieto-Quintic (siehe
\cite{BN}) entspricht.

Die Gruppe $\Gpq$ ist der Kern der Reduktionsabbildung
$\phi_q:\Gp\to\GL(4,\ZZ_q)$. Allgemein bezeichnen wir mit $\bar A$ die
Restklasse mod~$q$ irgendeiner Matrix~$A$.

\begin{lemma}\label{Sp4}
  Wenn $p\ne q$ ist, ist das Bild $\phi_q(\Gp)$ eine zu $\Sp(4,\ZZ_q)$
  konjugierte Untergruppe von~$\GL(4,\ZZ_q)$.
\end{lemma}
{\it Beweis:\/} F\"ur jedes $\bar\gamma\in\phi_q(\Gp)$ ist
${\bar{\gamma}}{\bar{\Lambda}}{\vphantom{\bar{\gamma}}^t{\bar{\gamma}}}
={\bar{\Lambda}}$. Da $p\ne q$ ist, ist $\bar\Lambda$ eine nicht-entartete
schief-symmetrische quadratische Form auf $\ZZ^4_q$ und damit zu
$J=\pmat{0&\bbEins_2\\ -\bbEins_2&0}$ \"aquivalent. Die Gruppe
$\Sp(\bar\Lambda, \ZZ_q)$ ist deshalb durch Konjugation mit $\bar{R}_p$ zu
$\Sp(4,\ZZ_q)$ isomorph und enth\"alt $\phi_q(\Gp)$.

Nach \cite[A.5.2]{F} ist $\Sp(4,\ZZ_q)$ von $J$ und $\pmat{\bbEins_2&S\\
  0&\bbEins_2}$, $S=\vphantom{S}^t{S}$, erzeugt. Also ist
$\Sp(\bar\Lambda,\ZZ_q)$ von $\bar R_pJ\bar R_p^{-1}=\pmat{0&\bar E^{-1}\\ 
  -\bar E & 0}$ und $\bar R_p\pmat{\bbEins_2&S\\ 
  0&\bbEins_2}\bar R_p^{-1}=\pmat{\bbEins_2&\bar{ S'}\\ 0&\bbEins_2}$
erzeugt, wobei $S'={\pmat{s_1 & s_2\\ s_2& s_3}} \cdot  E^{-1}$, $s_i\in\ZZ$.
Es seien ganze Zahlen $\lambda,\mu$ so gew\"ahlt, da{\ss} $\lambda p - \mu q =
1$ ist. Dann ist $\pmat{\bbEins_2 & \hat{S} \\ 0 &\bbEins_2}$ mit $\hat{S} =
\pmat{s_1 & \lambda s_2 \\ s_2 (1+\mu q) & \lambda s_3}$ in $\Gp$ und
insbesondere im Urbild von $\pmat{\bbEins_2& \bar{S'}\cr 0 & \bbEins_2}$.\\
W\"ahlen wir anderenfalls ganze Zahlen $\lambda,\mu$ mit $-\lambda p + \mu q = -1$ und
$\alpha,\beta$ mit $\alpha p + \beta q = -\mu$, so liegt
\[
M=\pmat{
0 & 0 & 1 & 0 \cr
0 & \beta q & 0 & \alpha q + \lambda \cr
-1 & 0 & 0 & 0 \cr
0 & -p & 0 & q}
\]
in $\Gp$ und  $\phi_q(M)=\pmat{0&\bar E^{-1}\cr -\bar E& 0}$.~\qed

\begin{lemma}\label{laGpqtransitiv}
  Ist $q\neq p$, so operiert $\phi_q(\Gp)$ transitiv auf
  $\ZZ_q^4\setminus\{0\}$. Ist $q=p$, so gibt es genau zwei
  $\phi_p(\Gp)$-Bahnen auf $\ZZ_p^4\setminus\{0\}$, und zwar
  $\ker\bar\Lambda\setminus\{0\}$ und $\ZZ_p^4\setminus\ker\bar\Lambda$.
\end{lemma}
{\it Beweis:\/} $\Sp(4,\ZZ_q)$ operiert transitiv auf
$\ZZ_q^{4}\setminus\{0\}$ und nach Lemma~\ref{Sp4} gilt dieses damit auch
f\"ur $\phi_q(\Gp)$, falls $q\neq p$ ist.

Im Falle $q=p$ ergibt sich: die Gruppe $\phi_p(\Gp)$ erh\"alt die Form
$\bar\Lambda$ sowie $\ker\bar\Lambda=\big\{(0, \bar v_2, 0, \bar
v_4)\mid \bar v_2, \bar v_4 \in \ZZ_p\big\}$. Die Untergruppe $\phi_p
j_2\big(\SL(2,\ZZ)\big)<\Gp$ operiert wie $\SL(2,\ZZ_p)$ auf
$\ker\bar\Lambda$, also transitiv.

Sei dann $(\bar v_1, \bar v_2, \bar v_3, \bar
v_4)\in\ZZ_p^4\setminus\{0\}$ und $(\bar v_1, \bar v_3)\neq
(0,0)$. Durch Anwendung eines geeigneten Elements aus $\phi_p
j_1\big(\SL(2,\ZZ)\big)$ k\"onnen wir annehmen, da\ss{} $(\bar v_1, \bar
v_3)= (1,0)$ ist. Aber $(1, \bar v_2, 0, \bar v_4)$ wird von
$\phi_q\left(M_2(-1)^{v_4}M_1(-1)^{v_2}\right)$, gefolgt von einem
Element aus $\phi_p j_1\big(\SL(2,\ZZ)\big)$ auf $v_0$
transformiert.~\qed

Ein Vektor $v\in\ZZ^4$ hei\ss{}t lang, wenn f\"ur alle $w\in\ZZ^4$ gilt:
$p|v\Lambda\vphantom{w}^t{w}$, sonst hei\ss{}t er kurz.
\begin{uremark}
Die langen Vektoren (bzw. kurzen Vektoren) sind gerade die Vektoren $v$ mit
$d_p(v) = p$ (bzw. $d_p(v) = 1$).
\end{uremark}
\begin{lemma}\label{kurzKlassen} Es seien $v, v'\in\ZZ^4$ kurz und
primitiv. Die Vektoren $v$ und $v'$ sind dann und nur dann
$\Gpq$-\"aquivalent, falls $\bar v=\overline{v'}$ ist.
\end{lemma}
{\it Beweis:\/} 
Falls $v'=\gamma v$ f\"ur ein $\gamma\in\Gpq$,
ist klar, da\ss{} dann $\bar v=\overline{v'}$ gilt.

Sei dann $\bar v=\overline{v'}$. Nach Lemma~\ref{laGpqtransitiv} d\"urfen
wir annehmen, da\ss{} $\bar v=v_0$ ist. Dann ist $\ggT(v_1,v_3)$ prim zu
$pq\ggT(v_2,v_4)$. Wenden wir Lemma~\ref{laggT} auf $(v_3,v_1,pq v_2,-pq
v_4)$ an, so liefert dieses ganze Zahlen $\lambda, \mu \in \ZZ$, so da\ss{}
$\ggT(v_1, v_3 + \lambda pq v_2 - \mu pq v_4) = 1$ ist. Wenden wir
$M_2(q)^{\mu}M_1(q)^{\lambda}$ auf $v$ an, so werden die Eintr\"age
$(v_1,v_3)$ nach $(v_1,v_3 + \lambda pqv_2 - \mu pqv_4 - \lambda\mu pq^2
v_1)$ transformiert. Es gilt $\ggT(v_1,v_3+\lambda pqv_2 - \mu pqv_4
-\lambda \mu pq^2 v_1) = 1$.

Wir k\"onnen also annehmen, da\ss{} $v_1$ und $v_3$ teilerfremd
sind. Da auch $\ggT(q,v_1)=1$ ist, gibt es ganze Zahlen $\lambda', \mu'$ so
da\ss{} $\lambda'v_1 + \mu'q v_3 = 1$ ist. Die Matrix
$\pmat{\lambda' & -v_3 \cr \mu'q & v_1}$ ist aus $\Gamma_1(q)$ und
$j_1\left(\pmat{\lambda' & -v_3 \cr \mu'q & v_1}\right)$ transformiert uns den
Vektor $v$ nach $(1,v_2,0,v_4)$. Anwendung von $M_2(q)^{-v_4/q}$
f\"uhrt diesen nach $(1,v_2,-pv_2v_4,0)$
\"uber. Wieder durch Anwendung eines geeigneten Elements aus
$j_1\big(\Gamma_1(q)\big)$ wird dieser nach $(1,v_2,0,0)$ transformiert, der
letztlich via $M_1(q)^{-v_2/q}$ nach $v_0$ \"ubergeht.~\qed

\begin{lemma}\label{langKlassen} Es seien $v, v'\in\ZZ^4$ lang und
primitiv. Ist $q\neq p$, so sind die Vektoren $v$ und $v'$ dann und nur dann
$\Gpq$-\"aquivalent, falls $\bar v=\overline{v'}$ ist. Ist $q=p$, so
zerf\"allt jede Restklasse mod~$q$ von langen Vektoren aus~$\ZZ^4$ in
genau $q^2$ unterschiedliche $\Gpq$-Bahnen.
\end{lemma}

Wie vorher k\"onnen wir nach Lemma~\ref{laGpqtransitiv} annehmen, da\ss{}
$\bar v=(0,1,0,0)\in\ZZ_q^4$. Da $\ggT(v_4,v_2,qv_1,qv_3)=1$ ist,
k\"onnen wir $\lambda, \mu\in\ZZ$ finden, so da\ss{}
$\ggT(v_2,v_4+\lambda qv_1+\mu qv_3)=1$. Durch Anwendung von zuerst
$M_2(q)^\lambda$ und danach $M_3(q)^\mu$ wird $v$ auf einen Vektor mit
$v_2, v_4$ teilerfremd transformiert. Durch Anwendung eines Elements
von $j_2\big(\Gamma_1(q)\big)$ k\"onnen wir auf $v=(v_1,1,v_3,0)$
reduzieren. Dann ist $p$ sowie $q$ Teiler von $\ggT(v_1,v_3)$, da $v$
lang, bzw. $\bar v=(0,1,0,0)$ ist.

Bei $q\neq p$ gilt also $pq|\ggT(v_1,v_3)$ und wir setzen
$v_1=pqv'_1$, $v_3=pqv'_3$. Wenden wir jetzt $M_2(q)^\nu$ gefolgt von
einem Element aus $j_2\big(\Gamma_1(q)\big)$ an, so k\"onnen wir statt
$(v_1,1,v_3,0)$ mit $(v_1,1,v_3+\nu pq,0)$ arbeiten. Wenn wir ein
geeignetes $\nu$ w\"ahlen, k\"onnen wir annehmen, da\ss{}
$v=(v_1,1,v_3,0)$ mit $\ggT(v'_1,v'_3)=1$ ist, und $v_3\neq 0$. Dieser
wird aber von $M_2(q)^{-v'_3}$ gefolgt von $M_3(q)^{v'_1}$ in
$(0,1,0,0)\in\ZZ^4$ transformiert.

Bei $q=p$ scheitert diese Argumentation, weil man statt $v_i=pqv'_i$
mit $\ggT(v'_1,v'_3)=1$ nur $v_1=pv'_1$, $v_3=pv'_3$ und
$\ggT(v'_1,v'_3)|p$ hat. Jeder lange primitive Vektor mit Restklasse
$\bar v=(0,1,0,0)\in\ZZ_p^4$ ist zu einem $(\alpha p,1,\beta p,0)$,
$0\le \alpha, \beta<p$ \"aquivalent. Es stellt sich noch die Frage, ob
diese Vektoren m\"oglicherweise unter sich $\Gpp$-\"aquivalent sind.
Der Vektor $(0,1,0,0)$ ist aber weder zu $(\alpha p,1,0,0)$ noch zu
$(0,1,\beta p,0)$ \"aquivalent ($\alpha,\beta\neq 0$). Sonst g\"abe es
Matrizen~$M^\pm$ aus $\Gpp$, deren zweite Zeile $(\alpha p,1,0,0)$,
bzw. $(0,1,\beta p,0)$ ist. Dann w\"are aufgrund der symplektischen
Bedingungen z.B. $M^+_{11}=M^+_{13}$: aber $\overline{M^+_{11}}=1$ und
$\overline{M^+_{13}}=0$. Die Bedingungen bei $M^-$ f\"uhren ebenfalls
zum Widerspruch. Da $(0,1,0,0)$ zu $(\alpha p,1,\beta p,0)$
$\Gp$-\"aquivalent ist und $\Gpp\triangleleft\Gp$ ist, gen\"ugt dieses
um zu beweisen, da\ss{} die Vektoren $(\alpha p,1,\beta p,0)$, $0\le
\alpha, \beta<p$ nicht $\Gpp$-\"aquivalent sind.~\qed

\begin{theorem}\label{LinieBahnen}
Bei $q\neq p,2$ (bzw. $q=p$, $q=2$) zerfallen die eindimensionalen
isotropen Unterr\"aume von $\bbQ^4$ unter der Operation von $\Gpq$ in
$q^4-1$ (bzw. $q^4-q^2$, $30$) Bahnen.
\end{theorem}
{\it Beweis:\/} 
Jede Gerade $\ell\subset\QQ^4$ enth\"alt genau zwei primitive Vektoren aus
$\ZZ^4$, die wir mit $\pm v_\ell$ bezeichnen. Zwei Geraden $\ell$ und
$\ell'$ sind genau dann \"aquivalent, wenn $\{\pm v_\ell\}$ mit $\{\pm
v_{\ell'}\}$ \"aquivalent ist, also wenn $v_\ell$ zu $\pm v_{\ell'}$
\"aquivalent ist. Bei $q\neq p$ hei\ss{}t das, da\ss{} $v_\ell$ und
$v_{\ell'}$ beide kurz oder beide lang sind, und
$\overline{v_\ell}=\overline{v_{\ell'}}$ ist. Die Bahnen der Geraden
ent\-sprechen also 
\[
[\overline{v_\ell},\mbox{ L\"ange}]\in(\ZZ_q^4\setminus\{0\})\times
\{\mbox{kurz, lang}\}/\pm 1.
\]
Bei $q\neq 2$ hei\ss{}t das $(q^4-1).2/2=q^4-1$ Bahnen. Bei $q=2$
operiert $\pm 1$ trivial und man hat also $(q^4-1).2=30$ Bahnen.

Bei $q=p$ und $v_\ell$ kurz l\"auft alles wie oben, ausgenommen nur,
da\ss{} $\overline{v_\ell}\in\ZZ_q^4\setminus\ker\bar\Lambda$ ist und es
daher f\"ur $\overline{v_\ell}$ statt $q^4-1$ nur $q^4-q^2$ m\"ogliche
Werte gibt. Man hat deshalb $(q^4-q^2)/2$ Bahnen von Geraden, die von einen
kurzen Vektor erzeugt sind. Ist $v_\ell$ lang, so h\"angt seine
$\Gpp$-\"Aquivalenzklasse nach Lemma~\ref{langKlassen} nicht nur von
$\overline{v_\ell}\in\ker\bar\Lambda\setminus\{0\}$ sondern auch von der
Restklasse mod~$p^2$ der Projektion von $v_\ell$ in
$\phi_p^{-1}(\ker\bar\Lambda)$ ab. Daf\"ur gibt es genau $p^2=q^2$
m\"ogliche Werte, man hat also $(q^2-1)q^2$ m\"ogliche
$\Gpp$-\"Aquivalenzklassen des Vektors $v_\ell$ und damit wiederum
$(q^4-q^2)/2$ Bahnen von Geraden, die von einem langen Vektor erzeugt sind,
also insgesamt $(q^4-q^2)$ Bahnen von Geraden.~\qed

Wenden wir uns nun den isotropen Ebenen aus $\bbQ^4$ zu. Nach
Korollar~\ref{corerzeugnis} kann jede isotrope Ebene $h$ als Erzeugnis
$v\wedge w$ eines kurzen Vektors $v$ und eines langen Vektors $w$ mit
$h_{\ZZ} = \ZZ v \oplus \ZZ w$, geschrieben werden. Notwendigerweise
ist $\overline{w}\neq \overline{v}$, da $h_\ZZ$ ebenfalls von $v-w$
und $w$ erzeugt ist und $v-w$ deshalb primitiv ist. W\"are
$\overline{w}=\overline{v}$, so w\"are $v-w$ durch~$q$ teilbar.

\begin{theorem} \label{LinienproEbene}
Ist $q\neq p,2$ (bzw. $q=p$, $q=2$), so enth\"alt jede isotrope Ebene bis
auf $\Gpq$-\"Aquivalenz genau $q^2-1$ (bzw. $q^2-q$, $6$) verschiedene
eindimensionale Unterr\"aume.
\end{theorem}
{\it Beweis:\/} Aufgrund der Transitivit\"at von $\Gp$ auf der Menge
der isotropen Ebenen reicht es aus, diese Aussage f\"ur die Ebene $h =
v_0 \wedge (0,1,0,0)$ zu zeigen. Diese Ebene enth\"alt kurze Vektoren
mit genau den Restklassen $(\alpha,\beta,0,0)$,
$(\bar\alpha,\bar\beta)\in \ZZ_q^2\setminus\{0\}$ und zwar
$(\alpha+\lambda q,\beta,0,0)$, $0<\alpha,\beta\le q$ mit $\alpha+\beta \neq 2q$ und
$\ggT(\alpha+\lambda q,\beta)=1$: bei $q=p$ aber ist $\bar\alpha=0$
nicht gestattet. Das hei\ss{}t bei $q\neq p,2$ (bzw.  $q=p$, $q=2$)
$q^2-1$ (bzw. $q^2-q$,~$3$) kurzen Vektoren bis auf
$\Gpq$-\"Aquivalenz, also $(q^2-1)/2$ (bzw. $(q^2-q)/2$,~$3$) von
kurzen Vektoren erzeugten Geraden. Bei $q\neq p$ (insbesondere bei
$q=2$) gilt dasselbe f\"ur lange Vektoren. Ist $q=p$, so sind die
\"Aquivalenzklassen von langen Vektoren auf~$h$ von $(\lambda
q^2+\alpha q, \beta, 0, 0)$ repr\"asentiert, wobei $0<\beta<q$, $0<
\alpha\le q$ und $\ggT(\lambda q^2+\alpha q, \beta)=1$ ist. Das
hei\ss{}t $q(q-1)$ langen Vektoren und wiederum $q(q-1)/2$ von langen
Vektoren erzeugten Geraden bis auf $\Gpp$-\"Aquivalenz.~\qed

\begin{lemma} \label{laInvEb}
  Es seien $h,h'$ zwei isotrope Ebenen. Gibt es zu jedem kurzen bzw. langen
  Vektor $v'$ aus $h'$ einen kurzen bzw. langen Vektor $v$ aus $h$, so
  da\ss{} $v$ und $v'$ $\Gpq$-\"aquivalent sind, so sind auch $h$ und $h'$
  $\Gpq$-\"aquivalent.
\end{lemma}
{\it Beweis:\/} Es reicht aus, die Behauptung f\"ur $h' = v_0 \wedge
(0,1,0,0)$ zu zeigen. Nach Voraussetzung enth\"alt $h$ einen kurzen Vektor
$v$ mit $\overline{v} = v_0$ und einen langen Vektor $w$ mit $\overline{w}
= (0,1,0,0)$ (und, falls $q=p$, auch $w_1\equiv w_3\equiv 0\mod q^2$). Ohne
Einschr\"ankung k\"onnen wir annehmen, da\ss{} das Gitter $h_\ZZ$ von $v$
und $w$ erzeugt wird. Es gibt nach Voraussetzung ein $\gamma \in \Gpq$, so
da\ss{} $v\cdot \gamma = v_0$ ist. Wir setzen $\tilde{w} :=  w \cdot \gamma$. Es folgt
aus der Isotropiebedingung $\tilde{w}_3 = 0$. Da $\tilde{h}_\ZZ=\ZZ v_0
v\oplus \ZZ \tilde{w}$ auch von $v_0$ und $\tilde{w} -\tilde w_1 v_0$
erzeugt wird, k\"onnen wir, ohne die Restklasse $\overline{\tilde w}$ oder
bei $q=p$ die Restklasse von $(\tilde w_1,\tilde w_3)\mod q^2$ zu \"andern,
$\tilde{w}_1 = 0$ voraussetzen. Wie $w$ ist $\tilde w$ primitiv, also
$\ggT(\tilde{w}_2,\tilde{w}_4)=1$. Durch Anwendung eines geeigneten
Elements aus $j_2\big(\Gamma_1(q)\big)$ wird $\tilde w$ auf $(0,1,0,0)$ und
$v_0$ wieder auf $v_0$, und damit $\tilde h=\tilde h_\ZZ\otimes\QQ$ auf
$h'$, gef\"uhrt.~\qed

\begin{theorem}\label{EbeneproLinie}
  Ist $q\neq 2$ (bzw. $q=2$), so liegt jeder eindimensionale Unterraum
  $\ell\subset\QQ^4$ in genau $(q^2-1)/2$ (bzw.~$3$) nicht
  $\Gpq$-\"aquivalenten isotropen Ebenen.
\end{theorem}
{\it Beweis:\/} Da $\Sp(\Lambda,\ZZ)$ transitiv auf den Mengen der kurzen
bzw.~langen Vektoren operiert, reicht es, die Aussage f\"ur den Fall $\ell
= \QQ v_0$ bzw.~$\ell = \QQ(0,1,0,0)$ zu beweisen. Es seien also $\ell =
\QQ v_0$ und $h = v_0 \wedge w$ eine isotrope Ebene, $w$ lang. Wie oben
k\"onnen wir ohne Einschr\"ankung $w=(0,w_2,0,w_4)$ mit $\ggT(w_2,w_4)=1$
voraussetzen. Dann k\"onnen wir verm\"oge $j_2\big(\Gamma_1(q)\big)$ auch
$0\le w_2,w_4 <q$ annehmen: nach Lemma \ref{langKlassen} sind diese Vektoren
(auch bei $q=p$) nicht \"aquivalent. Nach Lemma~\ref{laInvEb}
sind die entsprechenden Ebenen auch nicht \"aquivalent (aber $v_0\wedge
(0,w_2,0,w_4)=v_0\wedge (0,q-w_2,0,q-w_4)$). Das hei\ss{}t $(q^2-1)/2$
Ebenen, bzw. $3$~Ebenen bei $q=2$.

Die Aussage f\"ur den Fall $h= v \wedge (0,1,0,0)$ l\"a\ss{}t sich
(auch bei $q=p$) analog beweisen.

\begin{theorem}\label{AnzahlEbenen}
Es gibt bei $q\neq p,2$ (bzw. $q=p$, $q=2$) genau $(q^4-1)/2$
(bzw. $(q^2-1)(q^2+q)/2$,~$15$) Bahnen isotroper Ebenen.
\end{theorem}
{\it Beweis:\/} Das Titsgeb\"aude hat nach Satz \ref{LinieBahnen} und Satz
\ref{EbeneproLinie} $(q^2-1)(q^4-1)/2$ (bzw. $(q^2-1)(q^4-q^2)/2$,~$90$)
Kanten. Nach Satz \ref{LinienproEbene} entspricht jede Bahn isotroper
Ebenen $q^2-1$ (bzw. $q^2-q$,~$6$) Kanten.~\qed

Damit liegt f\"ur das Titsgeb\"aude ${\cal T}(\Gpq)$ eine
$(a_b,c_d)$-Konfiguration vor, wobei $a=q^4-1$ (bzw $q^4-q^2$,~$30$);
$b=(q^2-1)/2$ (bzw. $(q^2-1)/2$,~$3$); $c=(q^4-1)/2$ (bzw.
$(q^2-1)(q^2+q)/2$,~$15$); und $d=q^2-1$ (bzw. $q^2-q$,~$6$).

\begin{figure}[tbp]
  \begin{center}
\epsfbox{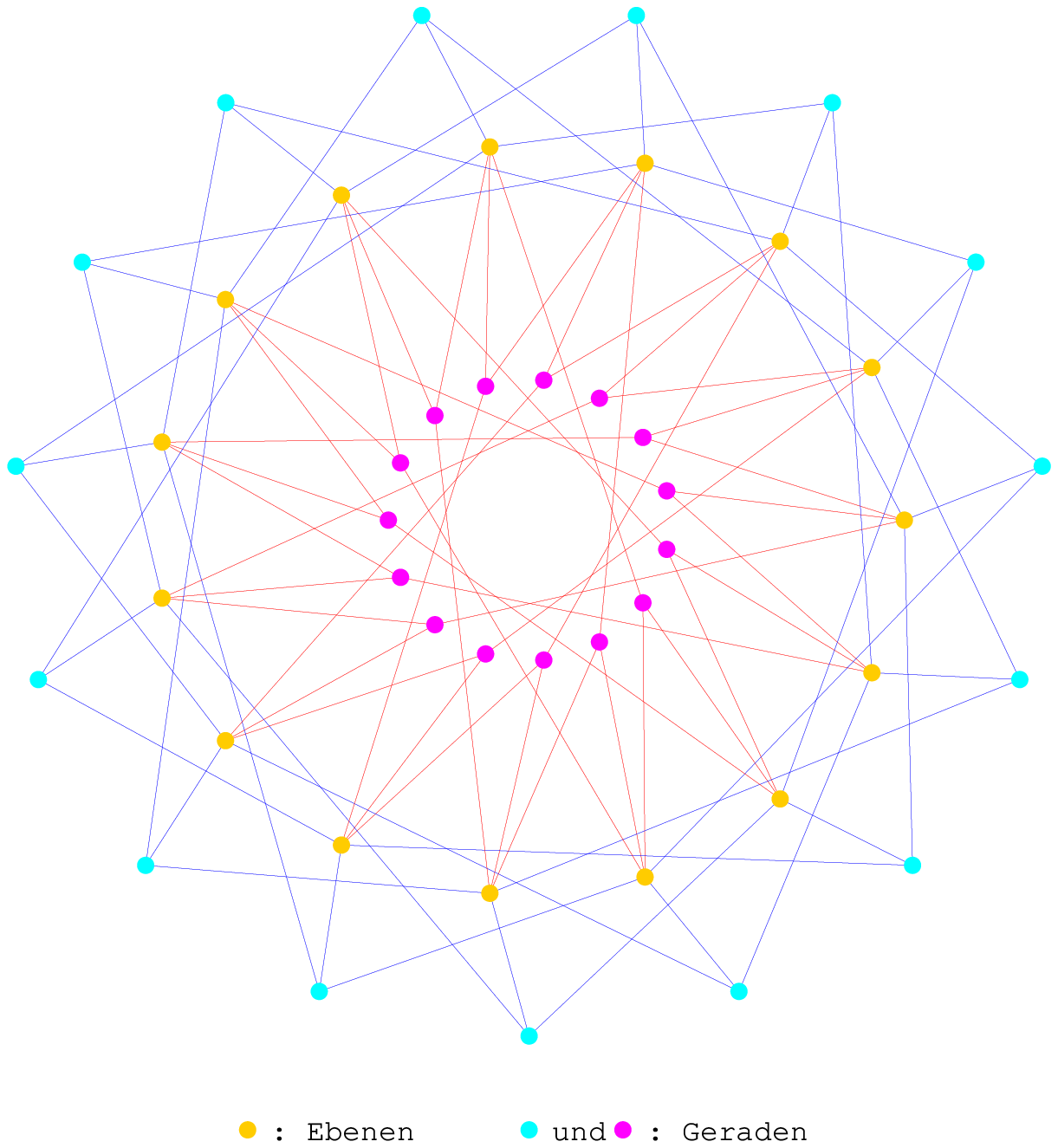} 
  \end{center}
\begin{center}
Bild 3: ${\cal T}({{\tilde{\Gamma}_{1,p}(2)}})$
\end{center}
\end{figure}
%
%$$
%BILD
%$$

Setzen wir nun $q=2$ oder $q=3$ und $q\neq p$ voraus. Wir haben gesehen,
da\ss{} je zwei eindimensionale Unterr\"aume genau dann \"aquivalent sind,
wenn die erzeugenden primitiven Vektoren entweder beide kurz oder beide
lang sind, sowie die gleiche Restklassen mod~$q$ bestimmen. Wir k\"onnen
also die Bahnen jeweils durch $\ZZ_q^4\setminus \{0 \}/\pm 1$ indizieren.
Ist $q=2$ oder $3$, so k\"onnen wir diese Menge auch als $\PP^3(\ZZ_q)$
auffassen. Nun k\"onnen wir jede isotrope Ebenen $h$ als Erzeugnis zweier
Vektoren $v$ und $w$ mit $\overline{v} \neq \overline{w}$ darstellen. Mit
der obigen Identifikation k\"onnen wir damit $\overline{h} = \overline{v}
\wedge \overline{w}$ als Gerade in $\Gr(1,\PP^3(\ZZ_q))$ auffassen. Nach
Lemma~\ref{laInvEb} sind zwei isotrope Ebenen $h$ und $h'$ genau dann
\"aquivalent, falls $\overline{h} = \overline{h'}$, das hei\ss{}t, wenn sie
die gleiche Gerade in $\Gr(1,\PP^3(\ZZ_q))$ bestimmen. Wir k\"onnen also
die Indexmenge der Bahnen isotroper Ebenen als Teilmenge der Grassmannschen
$\Gr(1,\PP^3(\ZZ_q))$ auffassen. Man kann jede Gerade in $\PP^{n}$ durch
ihre Pl\"uckerkoordinaten beschreiben. In unserem Fall sind diese durch
\[
 p_{ij} = \det \pmat{\bar v_i & \bar w_i \\
 \bar v_j & \bar w_j}
\]
f\"ur eine Gerade $\overline{v} \wedge \overline{w} \in
\Gr(1,\PP^3(\ZZ_q))$ gegeben: bekanntlich ist
$(p_{ij})$ genau dann in $\Gr(1,\PP^3(\ZZ_q))$, falls   
\[
 p_{12}p_{34} + p_{13}p_{24} + p_{14}p_{23} = 0.
\]
gilt.

Eine Gerade in $\Gr(1,\PP^3(\ZZ_2))$ repr\"asentiert genau dann
isotrope Ebenen, wenn $\overline{h} = \overline{v} \wedge
\overline{w}$ in $\ZZ_2$ isotrop ist, das hei\ss{}t, wenn die
Gleichung
\[
\det \pmat{\overline{v}_1 & \overline{w}_1 \cr
   \overline{v}_3 & \overline{w}_3 } = \det \pmat{\overline{v}_2 &
   \overline{w}_2 \cr
   \overline{v}_4 & \overline{w}_4 }
\]
erf\"ullt ist.

Eine Ebene ist also genau dann isotrop, wenn f\"ur ihre
Pl\"uckerkoordinaten $p_{13} = p_{34}$ in $\Gr(1,\PP^3(\ZZ_q))$ gilt. Die
Bahnen der isotropen Ebenen ent\-sprechen daher den L\"osungen in
$\PP^4(\ZZ_q)$ von $x_0x_1+x_2x_3+x_4^2=0$.

Im Fall $q=2$ kann man sich diese dann folgenderma\ss{}en veranschaulichen.
Den Raum $\PP^3(\ZZ_2)$ stellen wir uns durch die Punkte dar, die bei der
baryzentrischen Unterteilung eines Tetraeders entstehen. Dies sind genau
$4$~Eckpunkte, $6$~Kantenmittelpunkte, $4$~Fl\"achenmittelpunkte und ein
Schwer\-punkt, also ingesamt $15$~Punkte. Die Koordinaten der Eckpunkte seien
$(1,0,0,0)$, $(0,1,0,0)$, $(0,0,1,0)$ und $(0,0,0,1)$. Die Koordinaten der
anderen ergeben sich dann durch Addition. In dieses Gebilde k\"onnen wir
dann die 15 projektiven Geraden eintragen, die Bahnen isotroper Ebenen
entsprechen:

\begin{figure}[h]
  \begin{center}
\epsfbox{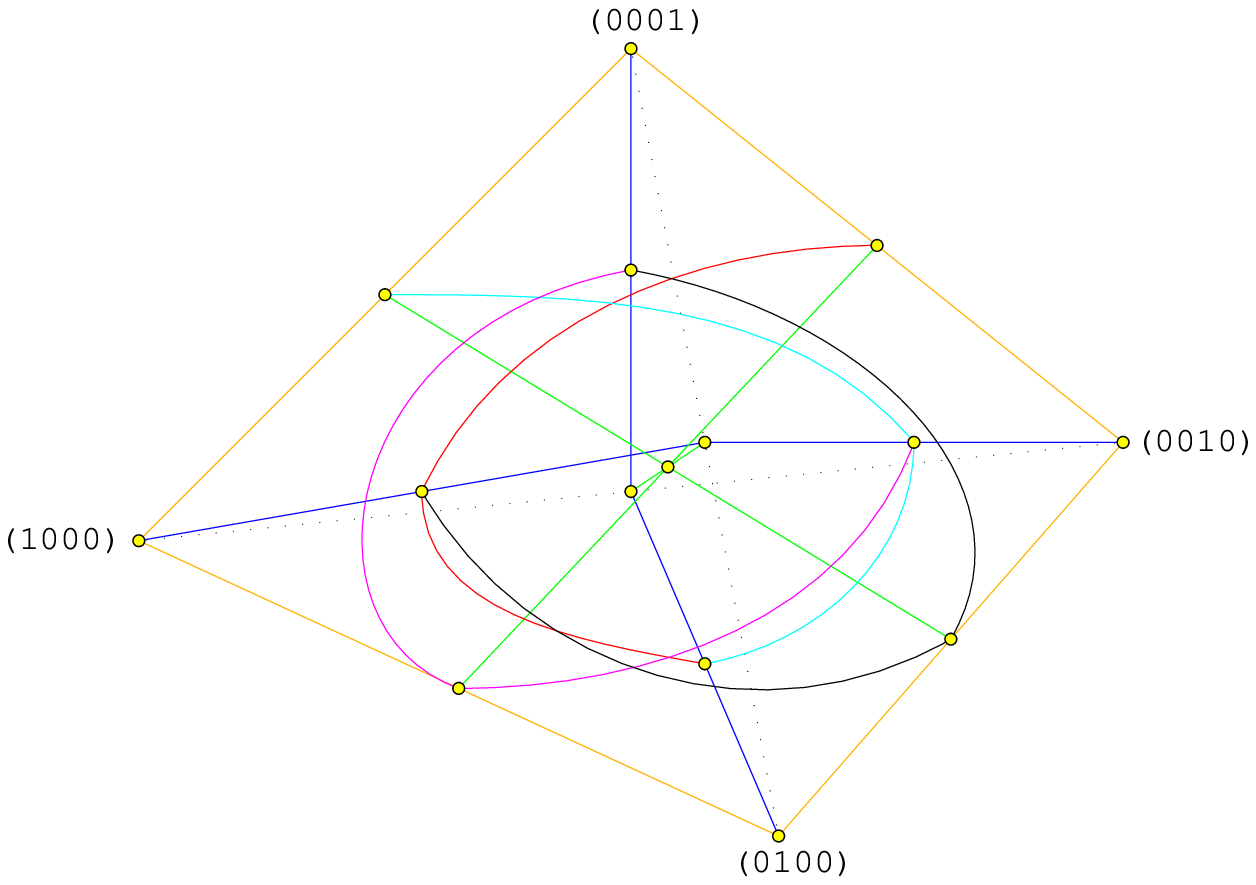} 
  \end{center}
\begin{center}
Bild 4: $\PP^3(\ZZ_2)$ und die Bahnen isotroper Ebenen im Fall $q=2$
\end{center}
\end{figure}
%
%$$
%BILD
%$$
%$$\epsfbox{../MetaPost/tetrahedron.1}$$

Anschrift der Autoren:

\bigskip

\begin{tabular}{ll}
M. Friedland & G.K. Sankaran\\
Institut f\"ur Mathematik & Department of\\
Universit\"at Hannover &  \ \ Mathematical Sciences\\
Postfach 6009 & University of Bath \\
D 30060 Hannover & Bath BA2 7AY\\
Germany & England\\
\\
{\tt friedland@math.uni-hannover.de} & {\tt gks@maths.bath.ac.uk}
\end{tabular}


\begin{thebibliography}{TITS}\label{biblio}
\frenchspacing 

\bibitem[BN]{BN} W.~Barth, I.~Nieto, \textit{Abelian surfaces of type
  $(1,3)$ and quartic surfaces with 16 skew lines}, J. Alg. Geom
  \textbf{3} (1994), 173--222.
  
\bibitem[F]{F} E.~Freitag, \textit{Siegelsche Modulfunktionen}, Grundlehren
  254, Springer, Berlin 1983.
  
\bibitem[Fr]{Fr} M.~Friedland, \textit{Das Titsgeb\"aude von Siegelschen
  Modulgruppen vom Geschlecht $2$}. Diplomarbeit, Hannover 1997.
  
\bibitem[H]{H} K.~Hulek, \textit{Elliptische Kurven, abelsche Fl\"achen und
  das Ikosaeder}, Jahresber. Deutsch. Math.-Verein. \textbf{91} (1989),
  126--147.

\bibitem[HW]{HW} K.~Hulek, S.~Weintraub, \textit{Bielliptic abelian surfaces},
Math. Ann. \textbf{283} (1989), 411--429.

\bibitem[HKW]{HKW} K.~Hulek, C.~Kahn \& S.~Weintraub, \textit{Moduli spaces
  of abelian surfaces: Compactification, degenerations and theta
  functions}. de~Gruyter 1993.
  
\bibitem[Sh]{Sh} G.~Shimura, \textit{Introduction to the arithmetic theory of
  automorphic function}, Princeton University Press 1971
  
\bibitem[Zi]{Zi} J.~Zintl, \textit{Invarianten von kompaktifizierten
  Mo\-dulr\"aumen polarisierter abelscher Fl\"achen}, Dissertation, Hannover

\end{thebibliography}
\end{document}